# THE RADIAL SPANNING TREE OF A POISSON POINT PROCESS[1]

By Francois Baccelli and Charles Bordenave

*Ecole Normale Supérieure and INRIA*

We analyze a class of spatial random spanning trees built on a realization of a homogeneous Poisson point process of the plane. This tree has a simple radial structure with the origin as its root.

We first use stochastic geometry arguments to analyze local functionals of the random tree such as the distribution of the length of the edges or the mean degree of the vertices. Far away from the origin, these local properties are shown to be close to those of a variant of the directed spanning tree introduced by Bhatt and Roy.

We then use the theory of continuous state space Markov chains to analyze some nonlocal properties of the tree, such as the shape and structure of its semi-infinite paths or the shape of the set of its vertices less than $k$ generations away from the origin.

This class of spanning trees has applications in many fields and, in particular, in communications.

**1. Introduction.** There is a current interest in graphs generated over a random point set. This class of graphs includes, in particular, the minimal spanning tree, the nearest neighbor graph or the geometric graphs [20, 30]. On these graphs, the theory has mainly focused on large sample asymptotics when the random point set consists of independent uniform points in the unit square.

In this paper we define a new model of rooted tree generated over a point set, the *radial spanning tree* (RST). This tree holds similarities with the minimal directed spanning tree [6, 21, 22] and the Poisson forest (as defined in [12, 13]). The mathematical analysis of the RST is conducted on a homogeneous Poisson point process (PPP) on the plane. The scope of this paper is two-fold: to prove some local geometric properties using stochastic

Received October 2005; revised September 2006.

[1]Supported in part by the EuroNGI network.

geometry and to derive global asymptotic properties. To this end, three embedded layers of analysis will be considered: local tree functionals around a vertex (e.g., the degree), properties of the path in the tree from a vertex to the root (e.g., the total length) or properties of the tree structure or topology (e.g., the characterization of the semi-infinite paths). These three layers of analysis rely on different probabilistic tools.

This class of spanning trees has applications in many fields and especially in communication networks. For instance, in wireless sensor networks, the sensed information is gathered at a central point called cluster head and flows from the sensors along the edges of a spanning tree, the root of which is the cluster head [10]. Similarly, in multicast communication, the broadcasting of information flows from the source to the members of the multicast group along the edges of a spanning tree rooted at the source [5]. In both cases, the nodes on which the spanning tree is built are randomly located in some domain (sensors may be randomly spread in a region and may then move in a random way; similarly, the members of the multicast group are a random subset of the set of users in a large network). In both cases, it is crucial that the decisions to build the spanning tree from the random set of nodes be local; in addition, the properties of the path between a leaf of the random spanning tree and its root determine the network performance. Another application with the same characteristics can be found within the context of a class of information location algorithms used in peer to peer networks [7].

In the next section we give the basic definition and we summarize our main results. Section 3 gives various distributions on local properties on the tree. In Section 4 we focus on the directed spanning forest (DSF) which can be seen as the limit of the RST far away from the origin. In Sections 5.1 and 5.3 we prove asymptotic shape theorems on the RST. Section 5.2 contains a proof of a law of large numbers on a semi-infinite path of the RST and Section 6 the proof of a law of large numbers for the spatial averages. Finally, in Section 7, we quickly discuss some extensions of the RST and give some open questions.

## 2. The radial spanning tree.

2.1. *Definition.* Let $|\bullet|$ denote a norm on $\mathbb{R}^d$ and $B(X, r)$ the open ball of radius $r$ and center $X$. A set of points $N$ of $\mathbb{R}^d$ is said to be *nonequidistant* if there do not exist points $X, Y, Z, T$ of $N$ such that $\{X, Y\} \neq \{Z, T\}$ and $|X - Y| = |Z - T|$.

If $N$ is a countable set of points in $\mathbb{R}^d$ with no accumulation points, we write, for all bounded sets $A$,

$$N(A) := \sum_{X \in N} \mathbb{1}(X \in A).$$

Let $N$ be a countable set of points in $\mathbb{R}^d$, nonequidistant, with no accumulation points and such that $O \in N$. We define the RST of $N$, $\mathcal{T} = (N, E)$, with $E$ the set of edges, as follows: each point, excluding the one in the origin, has an edge to its closest neighbor among the set of points which are closer to the origin. More formally, $\mathcal{T}$ is defined by

$$\text{If } |Y| < |X| \text{ and } X, Y \in N \text{ then } (X, Y) \in E$$

$$\iff \quad N(B(O, |X|) \cap B(X, |X - Y|)) = 0.$$

The nonequidistant property is needed to ensure that there is no tie: a vertex $X$ which is not the origin has exactly one nearest neighbor which is closer to the origin. We trivially deduce that $\mathcal{T}$ is a tree.

It is important to notice that the construction of the tree is local and that it does not minimize any global functional as do the minimal spanning tree or the trees analyzed by Howard and Newman in [17].

In this paper we will consider only point sets in the plane $\mathbb{R}^2$. All the results extend to higher dimension. If not otherwise specified, the norm $|\bullet|$ will be the Euclidean norm. We consider an orthonormal basis $(O, e_x, e_y)$.

Consider now some homogeneous Poisson point process $N$ on the plane, with intensity $\lambda > 0$ and in its Palm version: almost surely (a.s.) $O \in N$. Since the Poisson point process is a.s. nonequidistant, we can a.s. generate the RST $\mathcal{T}$ of $N$.

An instance of such RST is given in Figure 1.

Since $\mathcal{T}$ is scale-invariant, without loss of generality, we can set $\lambda = 1$, for a general $\lambda$, all results follow by multiplying distances by $\sqrt{\lambda}$. From the invariance of the PPP by rotation, we deduce also that the law of the RST is invariant by rotation.

We now define another random graph which will be used in the analysis of the asymptotic behavior of the RST. Let $(O, e_1, e_2)$ be an orthonormal basis of $\mathbb{R}^2$. On a nonequidistant point set $N$ with nonaccumulation, we define $\mathcal{T}_{e_1}$ the *directed spanning forest (DSF) with direction* $e_1$ as follows: the ancestor of $X \in N$ is the nearest point of $N$ which has a strictly larger $e_1$-coordinate. On a lattice, this random graph was first introduced by Bhatt and Roy in [6] and it also holds some similarities with the Poisson forest by Ferrari, Landim and Thorisson [12].

2.2. *Notation.* (A glossary of the main mathematical notation can be found at the end of the paper). The cardinality of set $S$ will also be denoted by $|S|$. For all points $X \in N \setminus O$, the ancestor of $X$ in the tree is the nearest neighbor which is closer to the origin. It will be denoted by $\mathcal{A}(X)$ and, by convention, we will take $\mathcal{A}(O) = O$. The iterate of order $k$ of $\mathcal{A}$ will be denoted by $\mathcal{A}^k(X)$.

We will denote by $R_O(X)$ the path from $X$ to the origin, defined as the sequence of ancestors of $X$. In Section 5.1 $R_O(X)$ will also be thought of as a piecewise linear curve in $\mathbb{R}^2$, namely, as the union of the edges connecting the points of this sequence.

Throughout the paper we will focus on functionals $F$ defined on the vertices of $\mathcal{T}$, such as the length of the edge $(X, \mathcal{A}(X))$, its orientation, and so on. For any $X \in \mathbb{R}^2$, $F(X)$ is defined on the tree $\mathcal{T}$ constructed over $N \cup \{X\}$. This definition is consistent with Slyvniak's theorem: if $P$ is the probability measure of the Poisson point process $N$, the Palm measure of $N$ with two points $O$ and $X$ is $P * \delta_X * \delta_O$. Hence, a.s. $N \cup \{X\}$ is nonequidistant and $F(X)$ can be interpreted as the value of the functional $F$ conditioned on the fact that $X$ is a vertex of the tree.

Several qualitative results of the present paper involve constants. For the sake of clarity, we will use $C_0$ to denote a positive constant to be thought of as small and $C_1$ to denote a positive constant to be thought of as large. The exact value of $C_0$ and $C_1$ may change from one line to the other and we could, for example, write $C_0/C_1 = C_0$. The important point is that $C_0$ and $C_1$ are universal constants that will never depend on the parameters of the problem.

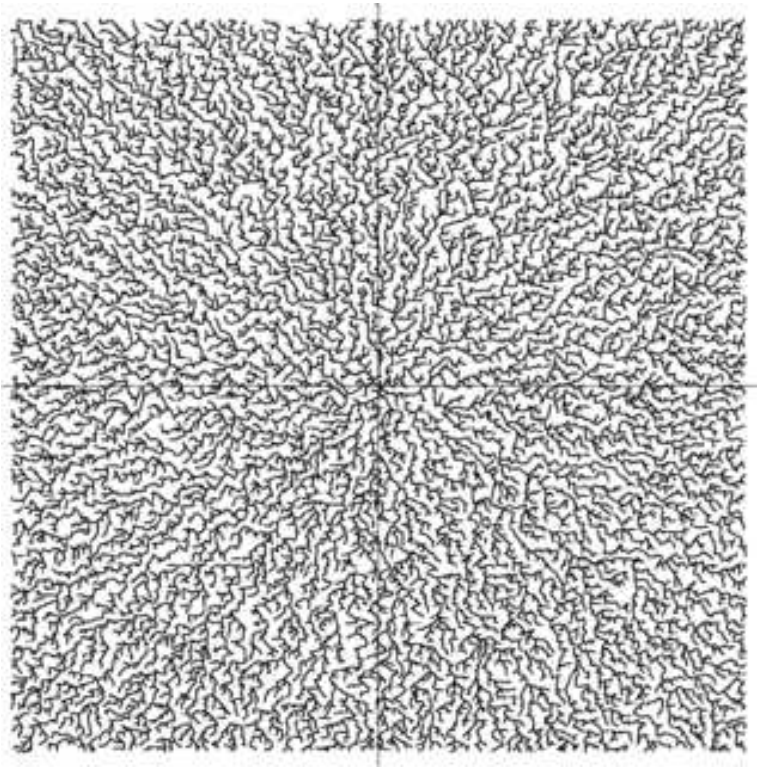

FIG. 1. *Radial spanning tree of* 25000 *points uniformly and independently distributed in the unit square.*

### 2.3. *Summary of results.*

2.3.1. *Topology of the radial spanning tree.* We have seen that if $N$ is a Poisson point process under its Palm version, $\mathcal{T}$ is a.s. a tree. We will check also that this tree is a.s. locally finite (i.e., no vertex has an infinite degree).

The next step to understand the intrinsic structure of $\mathcal{T}$ is to characterize its *ends*. An end is a semi-infinite self-avoiding path in $\mathcal{T}$, starting from the origin $(O = Y_0, Y_1, \ldots)$. The set of ends of a tree is the set of distinct semi-infinite, self-avoiding paths (two semi-infinite paths are nondistinct if they share an infinite sub-path).

A semi-infinite path $(0 = Y_0, Y_1, \ldots)$ has an asymptotic direction if $Y_n/|Y_n|$ has an a.s. limit in the unit sphere $S^1$.

The following theorem will be a consequence of Proposition 2.8 in [17]; it characterizes the ends of the RST:

THEOREM 2.1. *The following set of properties holds almost surely*:

- *every semi-infinite path has an asymptotic direction,*
- *for every $u \in S^1$, there exists at least one semi-infinite path with asymptotic direction $u$,*
- *the set of $u$'s of $S^1$ such that there is more than one semi-infinite path with asymptotic direction $u$ is dense in $S^1$.*

This theorem shows that the RST strongly differs from the minimal spanning tree. In dimension two it has been proved that the minimal spanning tree has only one end (see [2]).

Another property of the tree of interest to us is the set of points at tree-distance less than or equal to $k$ from the origin $\mathcal{T}(k) := \{X \in N : \mathcal{A}^k(X) = O\}$.

THEOREM 2.2. *There exists a constant $p > 0$, such that a.s., for all $\varepsilon > 0$, for $k$ large enough,*

$$(1) \qquad N \cap B(O, (1-\varepsilon)kp) \subset \mathcal{T}(k) \subset B(O, (1+\varepsilon)kp).$$

*Moreover, a.s. and in $L^1$, as $k$ goes to infinity,*

$$(2) \qquad \frac{|\mathcal{T}(k)|}{k^2} \to \pi p^2.$$

In other words, the graph-diameter of the RST generated by a Poisson point process inside a ball grows linearly with the number of points.

2.3.2. *Geometry of the radial spanning tree.* In this paper we will focus on two types of geometrical results.

The first type concerns functionals $F(X, \mathcal{T})$ which depend only on vertices around $X$. For example, we will give explicit formulae for the distributions of $L(X) := |X - \mathcal{A}(X)|$ (the length of the edge linking $X$ to its ancestor), $P(X) := |X| - |\mathcal{A}(X)|$ (the progress to the origin) and the expectations of the degree $D(X)$ of vertex $X$ in the tree.

More generally, we will prove a limit theorem for a large class of functionals, called *stabilizing functionals*. This class was first introduced by [18] and it was used by Penrose and Yukich [23, 24]; it is slightly modified here to suit to our framework. Roughly speaking, $F(X, \mathcal{T})$ stabilizes $\mathcal{T}$ if the value at $X$ depends only a small number of vertices around $X$.

Since $\mathcal{T}$ depends only on the point set $N$, with an abuse of notation, we can write equivalently $F(X, N)$ for $F(X, \mathcal{T})$ to make the dependence on $N$ explicit.

DEFINITION 2.3. Let $F(X, N)$ be a measurable $\mathbb{R}_+$-valued function defined on the vertices of a random graph $G = (N, E)$.

$F$ stabilizes $G$ if for all $X$ there exists an a.s. finite random variable $R(X) > 0$ such that $F(X, N \cap B(X, R(X)) \cup N') = F(X, N)$ for any locally finite point set $N' \subset \mathbb{R}^2 \setminus B(X, R(X))$ and the distribution of $R(X)$ does not depend on $X$.

We will prove the following theorem where $\mathcal{T}_{-e_x}$ denotes the DSF with direction $-e_x$:

THEOREM 2.4. *Let $F$ be a stabilizing functional for $\mathcal{T}_{-e_x}$. As $x$ tends to infinity, the distribution of $F(xe_x, \mathcal{T})$ converges in total variation toward the distribution of $F(O, \mathcal{T}_{-e_x})$.*

Note that the subscript $x$ in $e_x$ (the unit vector of the horizontal axis) has nothing to do with the real number $x$ that we let go to $\infty$.

This theorem has to be related to the convergence of geometric graphs as it is defined for the objective method (refer to [1]).

The second class of geometrical results is of a different nature: it concerns the path $R_O(X)$ from $X$ to the origin. The simplest result bears on $H(X) := \inf\{k : \mathcal{A}^k(X) = O\}$, the generation of $X$ in the RST, that is, the cardinal of $R_O(X)$. Let $X_0 = X, \ldots, X_{H(X)} = O$ be the sequence of points of $N$ in $R_O(X)$. Along the line of [17], it is interesting to look at $\sum_{k=0}^{H(X)-1} |X_{k+1} - X_k|^\alpha$, with $\alpha > 0$. More generally, we will prove the following:

THEOREM 2.5. *Let $p$ be the constant in Theorem 2.2. There exists a probability measure $\pi$ on $\mathbb{R}$ such that if $g(X)$ is a measurable function from $\mathbb{R}^2$ to $\mathbb{R}$, $|g(X)| \leq \max(C_1, |X|^\alpha)$ for some $\alpha > 0$, then a.s.*

$$\lim_{|X| \to +\infty} \frac{H(X)}{|X|} = \frac{1}{p} \quad and \quad \lim_{|X| \to +\infty} \frac{1}{H(X)} \sum_{k=1}^{H(X)} g([X_{k-1} - X_k]_{X_{k-1}}) = \pi(g),$$

*where $[U]_V$ is the vector $U$ rotated by an angle $-\theta$ and $V = r\cos\theta \cdot e_x + r\sin\theta \cdot e_y$ ($[U]_V$ is the vector $u$ expressed in the local coordinates of $V$).*

*If $g$ is continuous, we also have*

$$\lim_{|X| \to +\infty} \frac{1}{H(X)} \sum_{k=1}^{H(X)} g(X_{k-1} - X_k) = \pi(g), \qquad a.s.$$

We prove in Section 4 that the probability measure $\pi$ can be interpreted as the stationary measure on the infinite edge process in the DSF. Theorem 2.5 is a law of large numbers and $\pi$ can be understood as the limit probability measure of an edge conditioned on being on a semi-infinite path.

Theorems 2.4 and 2.5 are of a different nature. In particular, the mean value of $L(X)$ differs from the average of the lengths of the edges along the path $R_O(X)$. This can be explained by the fact that being on a long path is a bias. We will discuss this in Section 6.

## 3. Local properties of the radial spanning tree.

3.1. *Distribution of the length of edges.* Let $X \in \mathbb{R}^2$ and $L(X) := |X - \mathcal{A}(X)|$. Let $0 \leq r < |X|$; we get

$$(3) \qquad \begin{aligned} \mathbb{P}(L(X) \geq r) &= \mathbb{1}(r \leq |X|)\mathbb{P}(N(B(X, r) \cap B(O, |X|)) = 0) \\ &= \mathbb{1}(r \leq |X|)e^{-M(|X|, r)}, \end{aligned}$$

where $M(x, r)$ is the volume of the lens of the right part of Figure 2. Using the formula for the surface depicted by the left figure of Figure 2, we get that

$$(4) \qquad M(x, r) = x^2\left(\phi - \frac{\sin(2\phi)}{2}\right) + r^2\left(\frac{\pi}{2} - \frac{\phi}{2} - \frac{\sin(\phi)}{2}\right),$$

with

$$\phi = 2\arcsin\frac{r}{2x}.$$

With an abuse of notation, for $x > 0$, we define: $L(x) := L((x, 0))$. $L(x)$ has a density on $(0, x)$ equal to

$$\frac{d}{dr}M(x, r)e^{-M(x, r)}$$

and a mass at $x$ equal to

$$(5) \qquad e^{-M(x,x)} = e^{-x^2(2\pi/3 - \sin(2\pi/3))} = e^{-x^2(2\pi/3 - \sqrt{3}/2)}.$$

Notice that the distribution function of $L(x)$ is not stochastically monotone in $x$. Its mean $\mathbb{E}(L(x))$, which is plotted in Figure 3, is not monotone in $x$ either.

REMARK 3.1. In Appendix, Section A.1, we also consider the distribution of $L(X_n)$ when we number the points of $N$ by their a.s. increasing radius.

3.2. *Distribution of progress.* Given $L(X) = r < |X|$, consider the angle $\theta(X) := \widehat{OX\mathcal{A}(X)}$. Since $\psi = \pi/2 - \phi/2$ (see the right part of Figure 2), $\theta(X)$ is uniformly distributed on the interval $(\pi + \arg(X) - \psi, \pi + \arg(X) + \psi)$, with $\cos\psi = \sin(\phi/2) = r/(2|X|)$, that is, $\psi = \arccos\frac{r}{2|X|}$. Given $L(X) = |X|$, the angle $\theta(X)$ is $\pi + \arg(X)$.

The joint distribution density of $(L(X), \theta(X))$ is equal to

$$\mathbb{1}(r \in (0, |X|)) \frac{d}{dr} M(|X|, r) e^{-M(|X|, r)} \, dr$$

$$(6) \qquad \times \mathbb{1}(\theta \in (\pi + \arg(X) - \psi, \pi + \arg(X) + \psi)) \frac{d\theta}{2\psi}$$

$$+ \delta_{|X|}(r) \delta_{\pi + \arg(X)}(\theta) e^{-M(|X|, |X|)}.$$

The *progress* is defined as

$$P(X) := |X| - |\mathcal{A}(X)|.$$

The mean progress $P((x, 0))$ is plotted in function of $x$ in Figure 3.

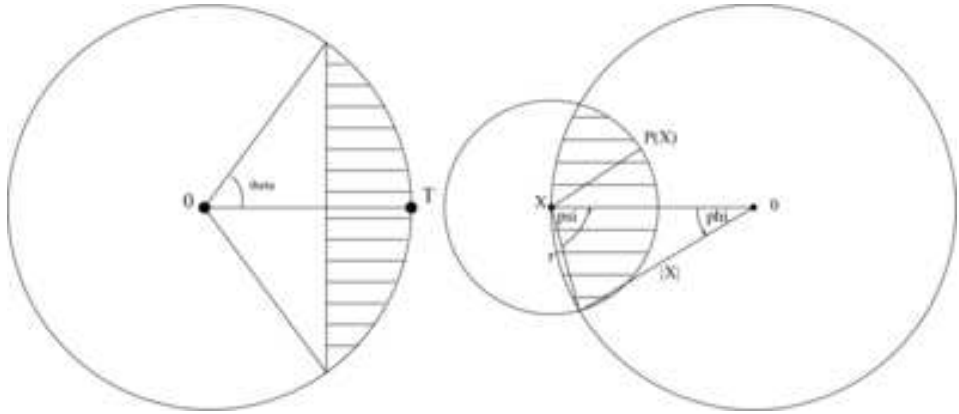

FIG. 2. Left: *the surface of the dashed lens is equal to* $\frac{|T|^2}{2}|2\theta - \sin 2\theta|$. Right: *the dashed lens.*

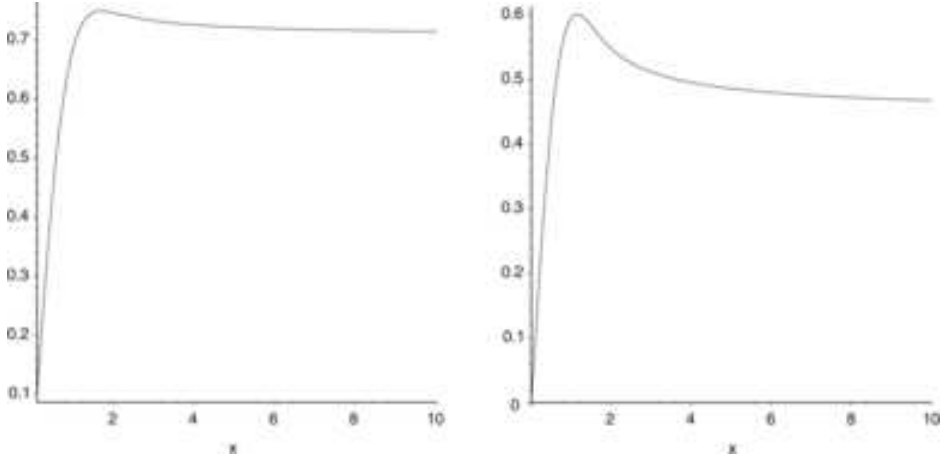

Fig. 3.  Left: *Mean of $L(x)$ in function of $x$. Right: Mean of $P(x)$ in function of $x$.*

### 3.3. *Mean degree of a vertex.*

#### 3.3.1. *Degree at the origin.* The degree of the origin is

$$D(O) = \sum_{T \in N \setminus \{O\}} \mathbb{1}(N(B(T, |T|) \cap B(O, |T|)) = 0).$$

Hence, using Campbell's formula, we get

$$(7) \qquad \mathbb{E}D(O) = 2\pi \int_0^\infty e^{-r^2(2\pi/3 - \sin(2\pi/3))} r \, dr = \frac{\pi}{2\pi/3 - \sqrt{3}/2} \approx 2.56.$$

The following property is also of interest. In the same vein as the boundedness of the degree of a vertex in the minimal spanning tree, we have the following:

LEMMA 3.2.  *The degree of vertex $O$ is bounded from above by $5$ a.s.*

PROOF.  Order the points directly attached to the origin by increasing the polar angle. Let $X$ and $Y$ denote two neighboring points in this sequence. Assume $|\vec{OX}| < |\vec{OY}|$. Denote by $\phi$ the angle between these two vectors. We have

$$|\vec{XY}|^2 = |\vec{OX}|^2 + |\vec{OY}|^2 - 2|\vec{OX}||\vec{OY}|\cos\phi.$$

Since $Y$ is attached to the origin, necessarily $|\vec{XY}|^2 > |\vec{OY}|^2$, which implies that

$$2|\vec{OX}||\vec{OY}|\cos\phi < |\vec{OX}|^2.$$

Using now the assumption that $|\vec{OX}| < |\vec{OY}|$, we get $\cos\phi < 1/2$. Hence, $|\phi| > \pi/3$.  $\square$

3.3.2. *Degree outside the origin.* The degree of a vertex $X \neq O$ is given by

$$D(X) = 1 + \sum_{T \in N \setminus \{X\}} \mathbb{1}(|T| \geq |X|) \mathbb{1}(N(B(T, |X - T|) \cap B(O, |T|)) = 0)$$

$$\times \mathbb{1}(O \notin B(T, |X - T|)). \tag{8}$$

Indeed, a point $T$ of modulus larger than $|X|$ shares an edge with $X$ if and only if there is no point of smaller modulus closer from $T$ than $X$.

Let $X \neq O$, $|X| = x$ and denoting $\mathbb{E}D(X)$ by $\mathbb{E}D(x)$, using Campbell's formula while taking the expectation of equation (8),

$$\mathbb{E}D(x) = 1 + \mathbb{E} \sum_{T \in N} \mathbb{1}(N(B(T, |X - T|) \cap B(O, |T|)) = 0)$$

$$\times \mathbb{1}(x \leq |T|) \mathbb{1}(|T| > |X - T|)$$

$$= 1 + \int_{\rho > x} \int_{-\arccos(x/(2\rho))}^{\arccos(x/(2\rho))} e^{-Q(x,\rho,\theta)} \rho \, d\rho \, d\theta,$$

where $Q(x, \rho, \theta)$ is the dashed surface in Figure 4 for $X = (x, 0)$ and $T = (\rho, \theta)$. The condition that $|T| > |X - T|$ [or, equivalently, that $\theta$ belongs to the interval $(-\arccos(\frac{x}{2\rho}), \arccos(\frac{x}{2\rho}))$] translates the fact that the origin should not be contained in this lens. Hence,

$$\mathbb{E}D(x) = 1 + \int_{\rho > x} \int_{-\arccos(x/(2\rho))}^{\arccos(x/(2\rho))} e^{-\rho^2/2|2\alpha - \sin 2\alpha|}$$

$$\times e^{-(\rho^2 + x^2 - 2\rho x \cos \theta)/2|2\beta - \sin 2\beta|} \rho \, d\rho \, d\theta,$$

where $\alpha$ and $\beta$ are the angles depicted in Figure 4.

If $u = \frac{\rho}{x}$, we have $\cos \alpha = (1 - u^{-2})/2 + u^{-1} \cos \theta$ and $\beta = (\pi - \alpha)/2$. Finally,

$$\mathbb{E}D(x) = 1 + 2x^2 \int_{u > 1} \int_0^{\arccos(1/(2u))} e^{-(u^2 x^2)/2(2\alpha - \sin 2\alpha)}$$

$$\times e^{-x^2/2(1 + u^2 - 2u \cos \theta)(\pi - \alpha - \sin \alpha)} u \, du \, d\theta. \tag{9}$$

The mean degree is plotted in Figure 5.

The following lemma is remarkable in view of Lemma 3.2.

LEMMA 3.3. *The degree of a vertex of the DSF is not bounded and in the RST a.s.*

$$\sup_{X \in N} D(X) = +\infty.$$

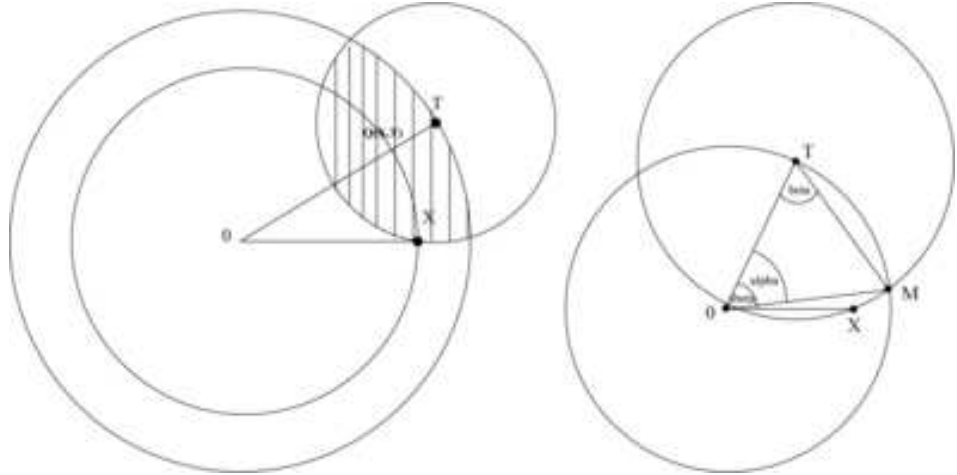

Fig. 4. Left: $Q(t, \rho, \theta)$. Right: *The $\alpha$ and $\beta$ angles.*

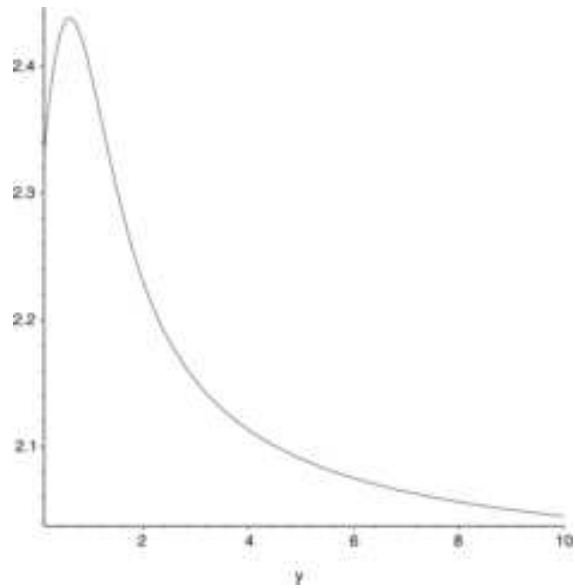

Fig. 5. $\mathbb{E}D(y)$ *as function of $y > 0$.*

PROOF. Let $\mathcal{A}_{-e_x}(X)$ be the ancestor of $X$ in the DSF with direction $-e_x$, $\mathcal{T}_{-e_x}$. The DSF built on the point set $\{X_n = (2^{-n}, 3^n), n \in \mathbb{N}\} \cup \{O\}$ gives, for all $n$, $\mathcal{A}_{-e_x}(X_n) = O$, in particular, the degree of the origin is infinite.

We now prove the second statement of the lemma.

Let $M \in \mathbb{N}^*$; for $n \geq 0$, we define $\mathbb{U}_n = [2^{-n} - \varepsilon, 2^{-n} + \varepsilon] \times [3^n - \varepsilon, 3^n + \varepsilon]$, $\mathbb{U}_{-1} = [-\varepsilon, \varepsilon] \times [-\varepsilon, \varepsilon]$, $A_M = B(O, 4^M) \setminus (\bigcup_{-1 \leq n \leq M} \mathbb{U}_n)$ and $E_M(X) = \{N(X + A_M) = 0, N(X + \mathbb{U}_n) = 1, -1 \leq n \leq M\}$. We have $\mathbb{P}(E_M(X)) = \delta > 0$ and if $|X - Y| > 2.4^M$, $E_M(X)$ and $E_M(Y)$ are independent (for $\varepsilon$ is small enough).

For $\varepsilon$ small enough, if $E_M(X)$ occurs, the point in $\mathbb{U}_{-1}(X)$ has degree at least $M$ in $\mathcal{T}_{-e_x}$. Similarly for the RST, if $|X|$ is large enough and if $E_M(X)$ occurs, the point in $\mathbb{U}_{-1}(X)$ has degree at least $M$ in $\mathcal{T}$.

Using the independence of the events $E_M(2k4^M e_x)$, $k \in \mathbb{N}$, we deduce that these events appear infinitely often and this concludes the proof.  □

3.4. *Limit distribution of an edge length.* Using equations (3) and (4), a direct computation gives

$$(10) \quad \lim_{|X| \to +\infty} \mathbb{P}(L(X) \geq r) = \exp\left(-\frac{\pi r^2}{2}\right),$$

$$\lim_{|X| \to +\infty} \mathbb{E}[L(X)] = \int_0^\infty e^{-\pi r^2/2}\, dr = \frac{1}{\sqrt{2}}.$$

By similar arguments, the asymptotic progress has for Laplace transform

$$(11) \quad \lim_{|X| \to +\infty} \mathbb{E}(e^{-sP(X)}) = \frac{1}{\pi} \int_{r=0}^\infty \int_{\theta=-\pi/2}^{\pi/2} e^{-sr\cos\theta} \exp\left(-\frac{\pi r^2}{2}\right) \pi r\, dr\, d\theta$$

$$= \int_{r=0}^\infty \int_{\theta=-\pi/2}^{\pi/2} e^{-sr\cos\theta} \exp\left(-\frac{\pi r^2}{2}\right) r\, dr\, d\theta.$$

In particular, the mean asymptotic progress is

$$(12) \quad \lim_{|X| \to +\infty} \mathbb{E}P(X) = \frac{1}{\pi} \int_{r=0}^\infty \int_{\theta=-\pi/2}^{\pi/2} \exp\left(-\frac{\pi r^2}{2}\right) \cos\theta\, dr\, d\theta = \frac{\sqrt{2}}{\pi}.$$

$L$ and $P$ are distributed as the length of the edge $(O, \mathcal{A}(O))$ and its $x$-coordinate in the DSF. This result is contained in Theorem 2.4 which we will prove in Section 3.6. With the same limit reasoning, for the degree, we get

$$\lim_{|X| \to +\infty} \mathbb{E}D(X) = 1 + \int_{\text{Half-plane}} \exp\left(-\frac{\pi |Y|^2}{2}\right) dY$$

$$= 1 + \int_{-\pi/2}^{\pi/2} \int_0^{+\infty} \exp\left(-\frac{\pi r^2}{2}\right) r\, dr\, d\theta$$

$$= 2.$$

3.5. *Expectation of the number of crossing edges.* Let $x > 0$ and $C(x)$ be the number of edges connecting a vertex belonging to the closure of $B(O, x)$ and a vertex outside. Almost surely, we can sort the points of $N$ by increasing norm $O = |T_0| < |T_1| < \cdots$. We have

$$C(|T_{n+1}|) = C(|T_n|) + D(T_{n+1}) - 2.$$

It follows that

$$C(x) = D(O) + \sum_{|T_k| \le x} (D(T_k) - 2).$$

Taking expectation, we deduce from Campbell's formula that

$$\mathbb{E}C(x) = \mathbb{E}D(O) + 2\pi \int_0^x (\mathbb{E}D(t) - 2)t\,dt.$$

In Section 3.4, we proved that $\lim_{x \to +\infty} \mathbb{E}D(x) = 2$. Hence, the evaluation of the asymptotic behavior of $\mathbb{E}C(x)$ requires that of $c = -\lim_{x \to +\infty} x^2 \mathbb{E}D(x)'$. We would then deduce $\lim_{x \to +\infty} \mathbb{E}C(x)/x = c$. The exact computation of $c$ is beyond the scope of this paper.

Let $\gamma$ be an arc on the circle of radius $r$ and center $O$ and let $C(r, \gamma)$ be the counting measure of the number of edges crossing $\gamma$. From the invariance by rotations of the PPP, we deduce that $\mathbb{E}C(r, \gamma) = l(\gamma)\frac{\mathbb{E}C(r)}{2\pi r}$, where $l(\gamma)$ is the length of the arc. In other words, the point process of edge crossings on the circle of radius $r$ is stationary and with intensity $\mu(r) = \frac{\mathbb{E}C(r)}{2\pi r}$.

3.6. *Proof of Theorem 2.4.* Let $X \in \mathbb{R}^2 \setminus \{O\}$ and $\mathcal{A}_{-e_1}(X)$ be the ancestor of $X$ in $\mathcal{T}_{-e_1}$ with direction $-e_1$ (a given vector). We define $\theta$ as the angle between $e_1$ and $X$. The next lemma compares $\mathcal{A}(X)$ with $\mathcal{A}_{-e_1}(X)$ if we build $\mathcal{T}$ and $\mathcal{T}_{-e_1}$ on the same PPP.

LEMMA 3.4. *Let $X \in \mathbb{R}^2 \setminus \{O\}$; there exists $C_1$ such that*

$$\mathbb{P}(\mathcal{A}(X) \ne \mathcal{A}_{-e_1}(X)) \le 1 \wedge C_1(1/|X| + \theta).$$

PROOF. Without loss of generality, we suppose $X = xe_x$, $x > 0$ and $\theta > 0$. The sets $L(X, e_1)$ and $K(X, e_1)$ are depicted in Figure 6. We note that if $\mathcal{A}_d(X) \notin K(X, e_1)$ and $\mathcal{A}(X) \notin L(X, e_1) \cup \{O\}$, then $\mathcal{A}_{-e_1}(X) = \mathcal{A}(X)$ and, hence,

$$\mathbb{P}(\mathcal{A}_{-e_1}(X) \ne \mathcal{A}(X)) \le \mathbb{P}(\mathcal{A}_{-e_1}(X) \in K(X, e_1))$$
$$+ \mathbb{P}(\mathcal{A}(X) \in L(X, e_1)) + \mathbb{P}(\mathcal{A}(X) = O).$$

The last term is easily computed: $\mathbb{P}(\mathcal{A}(X) = O) = e^{-M(x,x)}$.

To upper bound the second term, we notice that $L(X, e_1)$ is contained in a cone of angle $\theta$ (see Figure 6); hence, $\mathbb{P}(\mathcal{A}(X) \in L(X, e_1)) \le \theta/\pi$.

The first term is upper bounded similarly. Let $K_-(X, e_1)$ be the lower part of $K(X, e_1)$; we have

$$\mathbb{P}(\mathcal{A}_{-e_1}(X) \in K(X, e_1)) \leq 2\mathbb{P}(\mathcal{A}_{-e_1}(X) \in K_-(X, e_1))$$
$$\leq \int_0^\infty 2\pi r e^{-\pi r^2} \left( \arcsin\left(\frac{r}{2x}\right) + \theta \right) dr$$
$$\leq C_1(1/x + \theta). \qquad \square$$

We now prove Theorem 2.4. We recall that $F(X, \mathcal{T})$ is the value of the functional $F$ at vertex $X$ when the RST $\mathcal{T}$ is built on $N \cup \{X\}$. Without loss of generality, we can suppose that $X = xe_x$, $x > 0$. We build $\mathcal{T}$ and $\mathcal{T}_{-e_x}$ on the same PPP. For all $r > 0$, we define the event $J_r(X) = \{\mathcal{T} \cap B(X, r) = \mathcal{T}_{-e_x} \cap B(X, r)\}$. Let $F$ be a stabilizing functional for $\mathcal{T}_{-e_x}$. Using the terminology of Definition 2.3,

$$\mathbb{P}(F(X, \mathcal{T}) \neq F(X, \mathcal{T}_{-e_x}))$$
$$\leq \mathbb{P}(J_{R(X)}(X)^c)$$
$$\leq \mathbb{P}(R(X) > r) + \mathbb{P}(J_r(X)^c)$$
$$\leq \mathbb{P}(R > r) + \mathbb{P}\left( \bigcup_{T \in N \cap B(x,r)} \mathcal{A}(T) \neq \mathcal{A}_{-e_x}(T) \right)$$
$$\leq \mathbb{P}(R > r) + \mathbb{P}(N(B(x, r) \geq n) + nC_1(r/(x-r) + 1/(x-r))),$$

where we have used Lemma 3.4. For $\eta > 0$, we fix $r$ such that $\mathbb{P}(R > r) \leq \eta$. Note also that $\mathbb{P}(N(B(X, r)) > n) \leq \exp(-n \ln \frac{n}{e\pi r^2})$ (Lemma 11.1.1 of [28]).

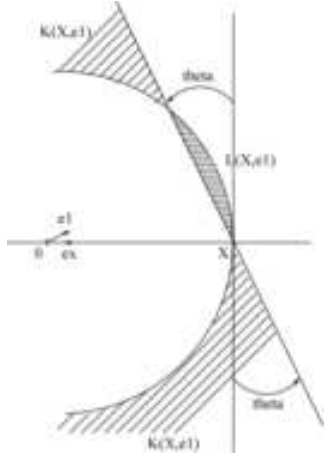

Fig. 6.  *The sets $L(X, e_1)$ and $K(X, e_1)$.*

Hence, taking $n = \lfloor x^\alpha \rfloor$ for some $0 < \alpha < 1$, we deduce that

$$\limsup_{x \to \infty} \mathbb{P}(F(X, \mathcal{T}) \neq F(X, \mathcal{T}_{-e_x})) \leq \eta,$$

and it follows $\lim_{|X| \to \infty} \mathbb{P}(F(X, \mathcal{T}) \neq F(X, \mathcal{T}_{-e_x})) = 0$.

To complete the proof, notice that $\mathcal{T}_d$ is stationary along $e_x$: $F(xe_x, \mathcal{T}_{-e_x})$ and $F(O, \mathcal{T}_{-e_x})$ have the same distribution.

REMARK 3.5. It is easy to check that the vector $(L(X), \theta(X))$ [and hence the progress $P(X)$], or the degree $D(X)$ are stabilizing functionals for $\mathcal{T}_{-e_x}$. So are the first $k$ segments of the path from $X$ to the origin in the DSF, for all finite $k$, or the subtree of the DSF rooted in $X$ and of depth $k$.

## 4. The directed path.
The directed spanning forest is the key tool for the asymptotic analysis of the RST. This section is dedicated to the analysis of the paths originating from some vertex.

4.1. *The underlying Markov chain.* We denote by $T_0, T_1, \ldots, T_n$ the sequence of the successive ancestors of point $T_0$ in the DSF with respect to the direction $-e_x$. For $i \geq 1$, the edge between $T_i$ and $T_{i-1}$ is $U_i := T_{i-1} - T_i$. Let $U_0 := -T_0$. From the definition of the tree, for $i \geq 1$ the $x$-coordinate of $U_i$ is positive (almost-surely) and we have:

$$T_n = -\sum_{l=0}^{n} U_l \quad \text{and} \quad T_n - T_k = -\sum_{l=k+1}^{n} U_l \qquad \text{for } k < n.$$

The conditional distribution of $U_{n+1}$ given $(U_0, \ldots, U_n)$ can be analytically determined. Indeed, let $\mathcal{H}$ be the half-plane $(x > 0)$, $D_0 := \mathcal{H}$ and for $n \geq 1$, let

$$\tag{13} D_n := \mathcal{H} \cap \left\{ \bigcup_{k=0}^{n-1} B(T_n - T_k, |U_{k+1}|) \right\}^c$$

(see Figure 9). Given $(U_0, \ldots, U_n)$, $U_{n+1}$ admits the following density in polar coordinates:

$$\tag{14} \frac{d}{dr} |B(O, r) \cap D_n| e^{-|B(O,r) \cap D_n|} \frac{\mathbb{1}_{(r,\theta) \in D_n} r \, dr \, d\theta}{\nu_1(C(O, r) \cap D_n)},$$

where $\nu_1$ denotes the 1-dimensional Lebesgue measure and $C(T, r)$ the circle of radius $r$ and center $T$.

Let $D_n$ be the set given by equation (13) for some fixed points $T_0, \ldots, T_n$. The following two lemmas will have important applications.

LEMMA 4.1.   *Let $(r_n, \theta_n)$ denote the coordinates of $U_n$. Then for all $0 < \alpha < \pi/2$,*

$$\mathbb{P}(|\theta_n| < \alpha | D_{n-1}, r_n) \geq \frac{\alpha}{\pi}.$$

*For $-\pi/2 \leq \alpha \leq \beta \leq \pi/2$, if the cone $\{(r, \theta); \theta \in [\alpha, \beta]\}$ is included in $D_{n-1}$, then*

$$\mathbb{P}(\theta_n \in [\alpha, \beta] | D_{n-1}, r_n) \geq \frac{\beta - \alpha}{\pi}.$$

PROOF.   Fix $r > 0$, in view of the geometry of $D_n$,

$$\left\{ \theta \in \left( -\frac{\pi}{2}, \frac{\pi}{2} \right) : (r, \theta) \in D_{n-1} \right\} = (-\theta_-(r), \theta_+(r)),$$

with $0 < \theta_\pm \leq \frac{\pi}{2}$.

From equation (14), for a fixed $r$, the p.d.f. of $\theta_1$ conditioned on $r_n = r$ is uniform on $(-\theta_-(r), \theta_+(r))$. If $\theta_- < \alpha$ and $\theta_+ < \alpha$, then $\mathbb{P}(|\theta_n| < \alpha | r_n = r, D_{n-1}) = 1$; else, supposing, for example, $\theta_+ \geq \alpha$, then

$$\mathbb{P}(|\theta_n| < \alpha | r_n = r, D_{n-1}) \geq \mathbb{P}(0 < \theta_n < \alpha | r_n = r, D_{n-1})$$
$$\geq \frac{\alpha}{\theta_+ + \theta_-} \geq \frac{\alpha}{\pi}.$$

The proof of the second assertion is similar.   □

Define the cones

$$c_\alpha = \begin{cases} X = (r, \theta) \in \mathbb{R}^2 : \theta \in [0, \alpha), & \text{for } \alpha > 0, \\ X = (r, \theta) \in \mathbb{R}^2 : \theta \in (\alpha, 0], & \text{for } \alpha < 0. \end{cases}$$

LEMMA 4.2.   *For all $n$,*

$$c_{\pi/6} \subset D_n \quad or \quad c_{-\pi/6} \subset D_n.$$

*In particular, if $(r_n, \theta_n)$ denote the coordinates of $U_n$, then*

$$\mathbb{P}(r_n \geq u | T_{n-1}, T_{n-2}, \dots, T_0) \leq e^{-\pi u^2/12}.$$

PROOF.   The proof relies on a simple geometrical argument. Suppose first $n = 2$ and consider a circle $C_2$ of radius $r_2$ and center $O$ (set to be $T_1$) and another circle $C_1$ of radius $r_1$ and center $T_0$, $O \in C_1$. In polar coordinates, $T_0$ is at $(r_1, \theta_1)$ and $T_2$ at $(\theta_2, r_2)$. For $|\theta_2 - \theta_1| \leq \frac{\pi}{2}$, the equation of $C_1$ in polar coordinates is $r = 2r_1 \cos(\theta - \theta_1)$ (see Figure 7). The point $T_2$ is somewhere on $C_2$; suppose, for example, that it is in the orthant $\theta_2 \in [-\frac{\pi}{2}, 0]$. If $\theta_1 \notin [-\pi, -\frac{\pi}{2}]$ or $r_1 \geq 2r_2 \cos \theta_1$, then $D_2$ contains $c_{\pi/2}$.

Suppose instead $\theta_1 \in [-\pi, -\frac{\pi}{2}]$ and $r_1 \le 2r_2 \cos\theta_1$ (see Figure 7). We have to prove that the largest cone with origin $T_2$ contained in $D_2$ contains $c_{\pi/6}$ or $c_{-\pi/6}$.

The worst case is when $T_2$ is at $M$, defined as the intersection of $C_1$ and $C_2$ in the orthant $\theta \in [-\frac{\pi}{2}, 0]$. We have $M = (r_2, \phi)$, with $\phi = \theta_1 + \arccos(\frac{r_2}{2r_1})$. In this case, an easy calculation shows that the largest cone contained in $D_2$ with origin $M = T_2$ is $\{X = M + (r, \theta) : \theta \in (-\frac{\pi}{2} - \theta_1 + 2\phi, \frac{\pi}{2} + \phi)\}$. The worst case is reached when $\phi = \frac{\theta_1}{3}$, since $\theta_1 \ge -\pi$. We deduce that $\max(\frac{\pi}{2} + \theta_1 - 2\phi, \frac{\pi}{2} + \phi) \ge \frac{\pi}{6}$.

This concludes the proof for $n = 2$. For $n \ge 3$, the largest cone contained in $D_n$ with origin $T_n$ is tangent to (at most) two circles, and the same conclusion holds. □

A sample path together with the associated exclusion discs are given in Figure 8.

From equation (14), the process $\{U_n\}, n \in \mathbb{N}$, is not Markov. We may circumvent this difficulty by defining $\tau_0 := 1$ and

$$(15) \quad \tau_{n+1} := \inf\left\{m > \tau_n : \mathcal{H} \cap \left\{\bigcup_{k=\tau_n-1}^{m-2} B(T_m - T_k, |U_{k+1}|)\right\}^c = \varnothing\right\}.$$

In words, $\tau_1$ is the first time $m$ such that no disc previous to $m-1$ intersects the half plane at $T_m$.

Each $\tau_n$ is a stopping time with respect to the process history generated by $\{U_m\}$. We call these times *Markovian times* (we will soon see why).

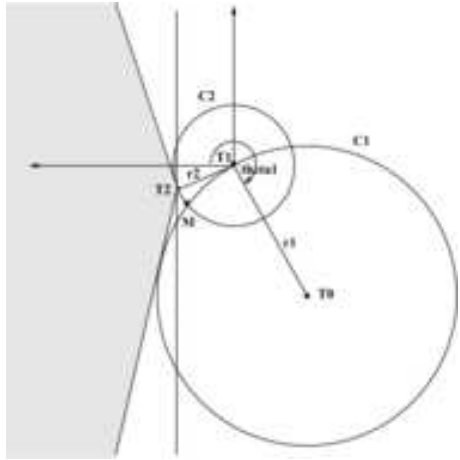

Fig. 7. *$D_2$ and the largest cone with origin $T_2$ contained in $D_2$.*

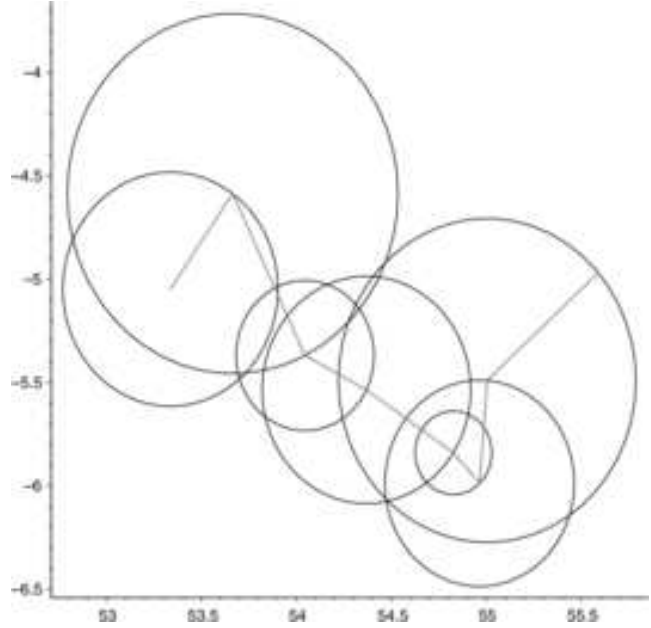

Fig. 8. *A sample path and the associated discs.*

Let $L_i := |U_{i+1}|$ denote the length of the edge from $T_i$ to $T_{i+1}$ and $P_i :=$ $\langle T_i - T_{i+1}, e_x \rangle$ the directed progress along this edge (where $\langle \cdot, \cdot \rangle$ denotes the usual scalar product).

The projection of $\mathcal{H} \setminus D_{\tau_0} = \mathcal{H} \cap \{B(T_1 - T_0, |U_1|)\}$ on the $x$-axis is an interval of length $\xi_1 = L_0 - P_0$. More generally, denote by $\xi_n$ the projection of $\mathcal{H} \setminus D_n = \mathcal{H} \cap \{\bigcup_{k=0}^{n-1} B(T_n - T_k, |U_{k+1}|)\}$ on the $x$-axis. This sequence satisfies the recurrence relation

$$(16) \qquad \xi_{n+1} = \max(\xi_n - P_n, L_n - P_n), \qquad n \geq 1,$$

so that, for all $n \geq 1$,

$$\xi_n = \max_{1 \leq i \leq n} \left( L_{i-1} - \sum_{k=i-1}^{n-1} P_k \right),$$

and the Markovian time $\tau_1$ is then simply rewritten as

$$\tau_1 = \inf\{m > 1 : P_m \geq \xi_m\}.$$

For instance, on the realization of Figure 9, $\tau_1 = 4$.

LEMMA 4.3. *For all $m \geq 1$,*

$$\mathbb{P}(\tau_m - \tau_{m-1} > n | D_{\tau_0}) \leq C_1 e^{-C_0 n} \ \text{and} \ \mathbb{P}(\tau_m > \beta m | D_{\tau_0}) \leq C_1 e^{-C_0 m},$$

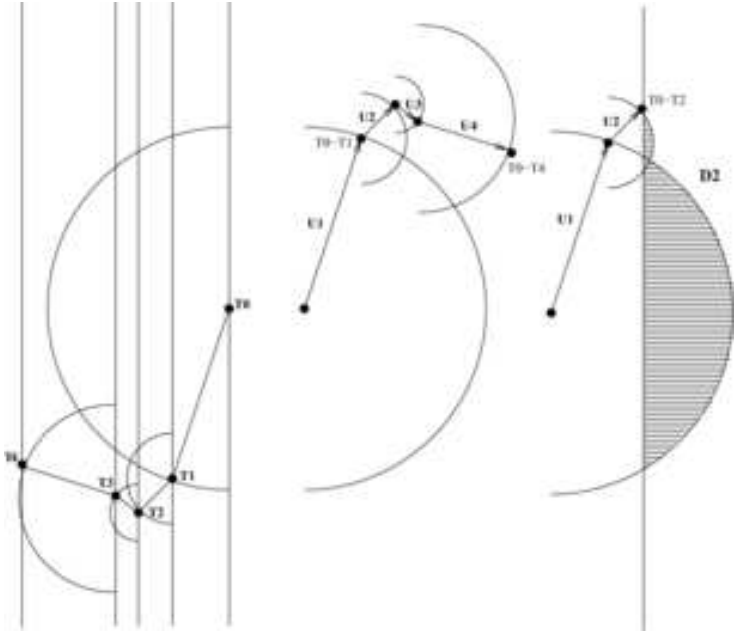

Fig. 9. *The process $T_0, T_1, \ldots, T_n$, the associated $U_1, \ldots, U_n$ and in dashed $\mathcal{H} \setminus D_2$.*

for some positive constants $\beta, C_0, C_1$ *(depending on $D_{\tau_0}$), so that, in particular, $\tau_m$ is a.s. finite.*

PROOF. This result follows from the statement of Lemma A.2 (in Appendix). We explain the connection between the setting of Lemma A.2 and the current setting in the particular case $m = 1$.

We first show that, for all $n$, the random variable $\xi_n$ defined in (16) is bounded from above by $\tilde{\xi}_n$, where $\tilde{\xi}_n$ is the maximal residual service time just before the $n$th arrival in a $GI/GI/\infty$ queue with i.i.d. service times $\{\tilde{L}_n\}$ and i.i.d. inter-arrival times $\{\tilde{P}_n\}$ where:

- $\tilde{L}_i$ is the distance from $T_i$ to the closest point in the cone (either $c_{\pi/6}$ and $c_{-\pi/6}$) that is fully included in $D_i$ if there is only one such cone. If both $c_{\pi/6}$ and $c_{-\pi/6}$ are included in $D_i$, then one selects one of them at random;
- $\tilde{P}_i$ is defined as follows: one samples an independent Poisson point process $N'$ of intensity 1 and one picks the point of

$$((N \cap D_i) \cup (N' \cap D_i^c \cap \mathcal{H})) \cap B(T_i, \tilde{L}_i)$$

which has the smallest progress from $T_i$.

Since for all $i$, $L_i \leq \tilde{L}_i$ and $P_i \geq \tilde{P}_i$ a.s., an immediate induction gives $\xi_i \leq \tilde{\xi}_i$ a.s., where $\tilde{\xi}_i$ is defined by the same recursion but with $\tilde{L}_i$ and $\tilde{P}_i$ in place of $L_i$ and $P_i$. This recursion is that of the maximal residual service time just before the arrivals in a $GI/GI/\infty$ queue.

Using a Loynes type argument (see, e.g., [4]), the sequence $\{\tilde{\xi}_i\}$ is easily seen to be stochastically monotone in $i$ and to converge weakly to a nondegenerate limit $\tilde{\xi}$ when $i$ tends to infinity.

Pick $a$ such that $\mathbb{P}(\tilde{\xi} \leq a) = \eta > 0$. Let $(r_n, \theta_n)$ be the polar coordinates of $U_n$. In view of the remark preceding the lemma,

$$\mathbb{P}(\tau_1 = m | D_0) \geq \mathbb{P}(\xi_m \leq a)\mathbb{P}(r_{m+1}\cos(\theta_{m+1}) \geq a | D_m)$$
$$\geq \mathbb{P}(\tilde{\xi} \leq a)\mathbb{P}(r_{m+1}\cos(\theta_{m+1}) \geq a | D_m).$$

Let $c(r, \alpha)$ denote the set $\{(\rho, \theta) : \rho < r; |\theta| < \alpha\}$. From Lemma 4.1, for all sets $D_m$, for all $\alpha \in (0, \frac{\pi}{2})$,

$$\mathbb{P}(r_{m+1}\cos\theta_{m+1} \geq a | D_m) \geq \mathbb{P}(r_{m+1}\cos\theta_{m+1} \geq a; |\theta_n| < \alpha | D_m)$$
$$\geq \frac{\alpha}{\pi}e^{-\nu_2(c(a/\cos\alpha, \alpha) \cap D_m)}$$
$$\geq \frac{\alpha}{\pi}e^{-a^2\alpha/\cos^2\alpha} = \mu > 0.$$

Let

$$\tilde{\tau}_1 = \inf\{m > 1 \, \tilde{P}_m \geq \tilde{\xi}_m\}.$$

Since $\tau_1 \leq \tilde{\tau}_1$, in order to prove the finiteness of $\tau_1$ (or the exponential bound on the tail of its law), it is enough to prove this property on $\tilde{\tau}_1$ and this is Lemma A.2. $\square$

For $n \geq 1$, we define the path between two successive Markovian times as the point process

$$(17) \qquad \Phi_n = \{0, T_{\tau_{n-1}} - T_{\tau_{n-1}+1}, \ldots, T_{\tau_{n-1}} - T_{\tau_n}\}$$

and we define $\Phi_0 = \{0, T_0 - T_1\}$. From Lemma 4.3, for all $n$, $\Phi_n$ is a finite point process on $\mathcal{X} = \mathbb{R}_+ \times \mathbb{R}$. We endow the set $\mathcal{N}_\mathcal{X}$ of finite point processes on $\mathcal{X}$ with the usual weak topology, $(\mathcal{N}_\mathcal{X}, \mathcal{B}(\mathcal{N}_\mathcal{X}))$, making a complete metric space (see [8] for details).

THEOREM 4.4. $\{\Phi_n\}, n \in \mathbb{N}$, *is a positive and uniformly ergodic Harris chain.*

The proof is organized in several steps. All definitions on Markov chains are taken from [19].

*Markov property.* Almost surely, we may write $\Phi_n = \{0, X_1^n, \ldots, X_{\tau_n}^n\}$, where the $x$-coordinates of the $X_i^n$ are increasing. By the definition of $\tau_n$, for all $n$, we have

$$D_{\tau_n} = \mathcal{H} \cap B(-U_{\tau_n}, |U_{\tau_n}|)^c.$$

It follows from equation [14] that the law of $U_{\tau_n+1} = X_1^n$ depends on $\Phi_0, \ldots, \Phi_n$ only through $U_{\tau_n} = T_{\tau_n-1} - T_{\tau_n}$.

More generally, since $\mathcal{H} \cap B(T_n - T_k, |U_{k+1}|) = \varnothing$ implies $\mathcal{H} \cap B(T_m - T_k, |U_{k+1}|) = \varnothing$ for $m \geq n$, it follows that, for all $m \geq \tau_n$,

$$D_m = \mathcal{H} \cap \left\{ \bigcup_{k=\tau_n-1}^{m-1} B(T_m - T_k, |U_{k+1}|) \right\}^c,$$

which depends only on $U_{\tau_n}, \ldots, U_m$. In particular, for all $n$,

$$\mathbb{P}(\Phi_{n+1} \in \bullet | \Phi_0, \ldots, \Phi_n) = \mathbb{P}(\Phi_{n+1} \in \bullet | U_{\tau_n}) = \mathbb{P}(\Phi_{n+1} \in \bullet | \Phi_n),$$

since $U_{\tau_n}$ is a functional of $\Phi_n$ only, which completes the proof of the Markov property.

In what follows, it will be useful to consider more general initial conditions than the two-point $\Phi_0$ defined above. Any finite point process $\Phi$ on $\mathcal{X}$ with $\Phi = \{0, X_1^0, \ldots, X_N^0\}$ and satisfying the constraint

$$D_N = \mathcal{H} \cap B(-U_N, |U_N|)^c$$

is an acceptable initial condition when taking $\tau_0 = N$. A special case is that where $\Phi = \{0\}$, where $N = 0$ and $D_0 = \mathcal{H}$.

In what follows, for a point pattern $\Phi$ as above, we will also use the notation $X_i(\Phi)$ for the $i$th point (the point, are ordered by their nondecreasing $x$-coordinates), $N(\Phi)$ for the number of points and $D(\Phi)$ for the set $D_N$.

As commonly done in Markov chain theory, we will denote by $P_\Phi$ the probability $\mathbb{P}(\bullet | \Phi_0 = \Phi)$ and by $P_\Phi^n$ the probability $\mathbb{P}(\Phi_n \in \bullet | \Phi_0 = \Phi)$, with $\Phi = \{0, X_1, \ldots, X_N\}$ a finite point process satisfying the above property. In particular, $P_0$ (or $P_0^n$) will denote the probability measure of the process conditioned on $\Phi_0 = \{0\}$.

Note that $P_\Phi$ depends on $\Phi$ only through $D(\Phi)$.

*Irreducibility and aperiodicity.* For $A \in \mathcal{B}(\mathcal{N}_\mathcal{X})$, define

(18) $$\sigma_A := \min\{n \geq 1 : \Phi_n \in A\}$$

and

$$L(\Phi, A) := P(\sigma_A < \infty | \Phi_0 = \Phi).$$

The irreducibility of a Markov chain on $\mathcal{N}_\mathcal{X}$ relies on the existence of a measure $\nu$ on $\mathcal{B}(\mathcal{N}_\mathcal{X})$ such that, for all $A$ in $\mathcal{B}(\mathcal{N}_\mathcal{X})$,

(19) $$\nu(A) > 0 \text{ implies } L(\Phi, A) > 0 \qquad \text{for all } \Phi \in \mathcal{N}_\mathcal{X}.$$

For all sub $\sigma$-algebras $\mathcal{G}$ of $\mathcal{B}(\mathcal{N}_\mathcal{X})$, let $P_\Phi^1|_\mathcal{G}$ denote the restriction of $P_\Phi^1$ to $\mathcal{G}$.

Let $\mathcal{F}$ denote the sub $\sigma$-algebra of $\mathcal{B}(\mathcal{N}_\mathcal{X})$ generated by the sets $\{N(\Phi) = 1\} \cap \{X_1(\Phi) \in B\}$, $B \in \mathcal{O}^+$, where $\mathcal{O}^+$ denotes the positive orthant $\mathbb{R}_+ \times \mathbb{R}_+$. We choose $\nu$ to be $P_0^1|_\mathcal{F}$.

As already pointed out, $P_\Phi^1$ depends on $\Phi$ only through $D(\Phi)$. In view of equation (14), if $D(\Phi) \cap \mathcal{O}^+ = \mathcal{O}^+$, then $P_\Phi^1|_\mathcal{F}$ and $P_0^1|_\mathcal{F}$ are equivalent measures. So if $D(\Phi) \cap \mathcal{O}^+ = \mathcal{O}^+$ and $P_0^1(A) > 0$, for some $A \in \mathcal{F}$, we then have $P_\Phi^1(A) > 0$, so that $L(\Phi, A) > 0$.

Consider now the case where $D(\Phi) \cap \mathcal{O}^+ \neq \mathcal{O}^+$. The first point of $\Phi_1$ is $X_1(\Phi_1) = U_{\tau_0+1}$. Let $(R_1, \Theta_1)$ be the coordinates of $X_1(\Phi_1)$. Let $(0, \xi)$ denote the projection on the $x$ axis of the set $\mathcal{H} \setminus D(\Phi)$. From equation (14), $\mathbb{P}(R_1 \cos \Theta_1 > \xi, \Theta_1 > 0) > 0$. If $R_1 \cos \Theta_1 > \xi$, $\tau_0 + 1$ is a Markovian time, so that $\tau_1 = \tau_0 + 1$, and $D(\Phi_1) = \mathcal{H} \cap B(-U_{\tau_1}, |U_{\tau_1}|)^c$. If, in addition, $\Theta_1 > 0$, $D(\Phi_1) \cap \mathcal{O}^+ = \mathcal{O}^+$, from what precedes, for $A \in \mathcal{F}$ and such that $P_0^1(A) > 0$, then $P_\Phi^2(A) > 0$, which implies $L(\Phi, A) > 0$.

The irreducibility is thus proved. The proof of the aperiodicity is along the same lines.

*Small set.* A set $S \in \mathcal{B}(\mathcal{N}_\mathcal{X})$ is small if there exists an integer $n > 0$ and a nontrivial measure $\nu$ on $\mathcal{B}(\mathcal{N}_\mathcal{X})$ such that, for all $A \in \mathcal{B}(\mathcal{N}_\mathcal{X})$,

$$(20) \qquad \inf_{\Phi \in S} P_\Phi^n(A) \geq \nu(A).$$

For $\Phi = \{0, X_1, \ldots, X_{\tau(\Phi)}\}$ in $\mathcal{N}_\mathcal{X}$, we define $(R(\Phi), \Theta(\Phi))$ as the coordinates of $X_\tau - X_{\tau-1}$. For $r > 0$ and $0 < \alpha < \pi/2$, let $\mathcal{S}(r, \alpha) \in \mathcal{B}(\mathcal{N}_\mathcal{X})$ be defined by

$$(21) \qquad \mathcal{S}(r, \alpha) = \{\Phi \in \mathcal{N}_\mathcal{X} : |\Theta(\Phi)| \leq \alpha, R(\Phi) \leq r\}.$$

Let $\mathcal{G}$ denote the sub $\sigma$-algebra of $\mathcal{B}(\mathcal{N}_\mathcal{X})$ generated by the sets $\{\tau(\Phi) = 1\} \cap \{X_1(\Phi) \in B\}$, $B \in C_\alpha$, where $C_\alpha$ denotes the cone

$$C_\alpha = \{(r, \theta) \in \mathbb{R}^2 : \theta \in [-\pi/2 + \alpha, \pi/2 - \alpha]\}.$$

We now prove that

$$(22) \qquad \mathcal{S} = \mathcal{S}(r, \alpha)$$

is a small set for $\nu = g(r, \alpha) P_0^1|_\mathcal{G}$ and $n = 1$, where $g(r, \alpha)$ is a positive constant to be determined below.

From equation (14), for all $\Phi \in \mathcal{S}$,

$$\begin{aligned}
P_\Phi(X_1(\Phi_1) \in C_\alpha, \tau(\Phi_1) = 1) &= P_\Phi(|\Theta_1| \leq \pi/2 - \alpha, R_1 \cos \Theta_1 \geq r) \\
&\geq P_0(|\Theta_1| \leq \pi/2 - \alpha, R_1 \cos \Theta_1 \geq r) \\
&\geq \frac{\pi/2 - \alpha}{\pi} e^{-r^2(\pi/2-\alpha)/\sin^2 \alpha} = g(r, \alpha) > 0.
\end{aligned}$$

Note also that from equation (14), for all $B \subset \mathcal{E}$,

$$P_\Phi(X_1 \in B \mid |\Theta_1| \leq \pi/2 - \alpha, R_1 \cos \Theta_1 \geq r)$$
$$= P_0(X_1(\Phi_1) \in B \mid |\Theta_1| \leq \pi/2 - \alpha, R_1 \cos \Theta_1 \geq r).$$

These two remarks imply that, for all $A \in \mathcal{G}$,

$$P_\Phi(\Phi_1 \in A) \geq g(r, \alpha) P_0(\Phi_1 \in A),$$

which concludes the proof.

*Positivity and uniform ergodicity.* Let $\sigma_\mathcal{S} = \inf\{n \geq 1 : \Phi_n \in \mathcal{S}\}$ be the first return time to $\mathcal{S}$ with $\mathcal{S}$ defined as above.

LEMMA 4.5. *For all $r > 0$, $\mathcal{S} = \mathcal{S}(r, \pi/6)$, we have $\sup_\Phi E_\Phi \sigma_\mathcal{S} < \infty$. More precisely, $\sup_\Phi \mathbb{P}(\sigma_\mathcal{S} > n) \leq (1 - \delta(r))^n$, with $\delta(r) = (1 - e^{-\pi r^2/12})/6$.*

PROOF. From Lemmas 4.1 and 4.2, for all $n$,

$$\mathbb{P}\left(r_n \leq r, |\theta_n| \leq \frac{\pi}{6} \Big| T_{n-1}, \ldots, T_0\right) \geq \frac{1}{6}(1 - e^{-\pi r^2/12}).$$

Hence, for all $\Phi$,

(23) $$P_\Phi\left(R(\Phi_1) \leq r, |\Theta(\Phi_1)| \leq \frac{\pi}{6}\right) \geq \frac{1}{6}(1 - e^{-\pi r^2/12}),$$

from which one easily deduces that $\sup_\Phi E_\Phi \sigma_\mathcal{S} < \infty$. $\square$

In view of Theorems 10.4.10 and 16.0.2 of [19], $\{\Phi_n\}$ is positive and uniformly ergodic.

4.2. *Limit theorems.* We follow the approach of Athreya and Ney (see [3]). Since the Markov chain $\Phi$ is Harris recurrent and strongly aperiodic, we can build an increasing sequence of finite stopping times $N_k$ (on an enlarged probability space) such that $N_0 = 0$ and

(24) $$P_{\Phi_0}(\Phi_n \in A, N_k = n) = \mu(A) P_{\Phi_0}(N_k = n),$$

where $\mu$ is a probability measure on $\mathcal{N}_\mathcal{X}$. Therefore, $N_k - 1$ is a regenerative time: the sequences $(\Phi_{n+N_k}), n \in \mathbb{N}$, and $(\Phi_n), 0 \leq n \leq N_k - 1$, are independent.

Lemma 4.5 implies that $N_k - N_{k-1}$ is stochastically dominated by a geometric law. We define $\theta_k = \tau_{N_k}$ with $\tau_n$ defined in (17); since $\xi_{\theta_k} \leq r$ [with $\xi_n$ defined in (16)], from Lemma 4.3, $\mathbb{P}(\theta_{k+1} - \theta_k \geq t | \mathcal{F}_{\theta_k}) \leq C_1 \exp(-C_0 k)$, for some positive constants.

Finally we have the following:

THEOREM 4.6. *There exists an increasing sequence of finite stopping times $\{\theta_k\}, k \in \mathbb{N}$, (on an enlarged probability space) such that $\theta_0 = 0$, for all $k \geq 0$, $\mathbb{P}(\theta_{k+1} - \theta_k > n | \mathcal{F}_{\theta_k}) \leq C_1 \exp(-C_0 n)$, for some positive constants and the vectors $(U_{\theta_k+1}, \ldots, U_{\theta_{k+1}}), k \in \mathbb{N}^*$ are i.i.d.*

For $\Phi = \{0, X_1, \ldots, X_\tau\}$ in $\mathcal{N}_{\mathcal{X}}$, we define $f(\Phi) = \sum_{n=1}^{\tau} \max(1, |X_n|^\alpha)$, $\alpha > 1$. Using the upper bounds of Lemma 4.3 and Lemma 4.2, we get

$$(25) \qquad \sup_{X \in \mathcal{S}} \mathbb{E}_X \left( \sum_{k=0}^{\sigma_{\mathcal{S}}-1} f(\Phi_k) \right) < \infty,$$

where $\sigma_{\mathcal{S}}$ is defined by equation (18) and $\mathcal{S}$ is the small set defined above.

Let $\Pi$ be the invariant distribution of the Markov chain. We define the *invariant distribution* of the edge process $T_0 - T_1, \ldots, T_{n-1} - T_n, \ldots$ as

$$(26) \qquad \pi(A) = \mathbb{E}_\Pi(\tau)^{-1} \mathbb{E}_\Pi \left( \sum_{l=0}^{\tau-1} \mathbb{1}((X_{l+1} - X_l) \in A) \right), \qquad A \in \mathbb{R}^2.$$

We may now deduce a first limit theorem.

THEOREM 4.7. *For all measurable functions $g$, let $S_n = \sum_{k=0}^{n-1} g(T_{k+1} - T_k)$. If $g(X) \leq \max(C, |X|^\alpha)$ for $C > 0$, $\alpha > 0$, then a.s.,*

$$\lim_{n \to \infty} \frac{S_n}{n} = \mathbb{E}_\Pi(\tau)^{-1} \mathbb{E}_\Pi \left( \sum_{l=0}^{\tau-1} g(X_{l+1} - X_l) \right) = \pi(g).$$

*Let $G(\Phi) = \sum_{l=0}^{\tau-1} (g(X_{l+1} - X_l) - \pi(g))$; if*

$$\gamma^2 = \mathbb{E}_\Pi(G(\Phi_0)^2) + 2 \sum_{k=1}^{\infty} \mathbb{E}_\Pi(G(\Phi_0)G(\Phi_k)) > 0,$$

*then a central limit theorem also holds:*

$$\frac{1}{\gamma\sqrt{n}}(S_n - \pi(g)) \xrightarrow{d} \mathcal{N}(0, 1).$$

*If $\gamma = 0$, then $\frac{1}{\sqrt{n}}(S_n - \pi(g))$ tends a.s. to 0.*

This theorem characterizes the asymptotic behavior of a path in the directed spanning forest.

PROOF. Theorem 4.7 is a direct application of Theorem 17.0.1 of [19]. Theorems 14.2.3, 14.2.4 and 14.3.7 of [19] along with equation (25) ensure that all requirements are fulfilled. □

COROLLARY 4.8. *There exists positive constants $p$, $p_y$ and for $\alpha > 0$, $l_\alpha$ such that*

$$\lim_{n\to\infty} \frac{1}{n} \sum_{k=0}^{n-1} \langle T_k - T_{k+1}, e_x \rangle = p_x = p, \tag{27}$$

$$\lim_{n\to\infty} \frac{1}{n} \sum_{k=0}^{n-1} |\langle T_k - T_{k+1}, e_y \rangle| = p_y, \tag{28}$$

$$\lim_{n\to\infty} \frac{1}{n} \sum_{k=0}^{n-1} |T_{k+1} - T_k|^\alpha = l_\alpha. \tag{29}$$

By simulation of 20000 transitions of the chain, one obtains that $p \sim 0.504$, $p_y \sim 0.46$ and $l_1 \sim 0.75$. The value of $p$ is significantly larger than the mean asymptotic progress as evaluated in (12). The latter can be seen as the expectation of the progress from point $T_0 = 0$ under the Palm probability of the Poisson point process $N$, whereas $p$ is the expectation of the progress under another measure: for all $n$, let $P_n = |U_n \cdot e_x|$ be the progress from the $n$th ancestor of point $T_0 = 0$. From the above analysis, the law of $P_n$ under the Palm probability of $N$ converges weakly to a limit when $n$ tends to infinity, and $p$ is the mean of the limit law. A similar observations holds for $l_1$ when compared to (10): as it is the case for progress, on a long path, the magnitude of the hop from a point to its ancestor is "boosted" by the presence of its offspring.

Note that it is also possible to derive a functional central limit theorem for the sequence $(T_k)_k$ from Theorem 17.4.4 of [19].

REMARK 4.9. An interesting question along the line of [12] and [15] is whether the directed spanning forest is almost surely a tree, namely, whether two sequences of ancestors coalesce almost surely (in dimension 2 and 3). Similarly, every vertex has almost surely a finite number of successors. One can also expect that the directed spanning tree converges toward the Brownian web (see [13, 14, 29]). If this holds true, then the directed spanning tree has only one semi-infinite path, whereas in view of Theorem 2.1, the RST has a semi-infinite path in every direction. The edge process of successive ancestors in the RST converges in some sense toward the edge process of successive ancestors in the directed spanning stree; however, the two trees have a completely different topology.

4.3. *Maximal deviation.* We end this section with a result on the deviation of the path from its mean.

Let $R(x)$ or $R(X)$ denote the path from $X = (x, 0)$ in the DSF with direction $-e_x$; $R(x)$ may be parametrized as a piecewise linear curve $(t, Y(t))_{t \le x}$ in $\mathbb{R}^2$. The maximal deviation of this curve between $x'$ and $x$ with $x' \le x$ is

defined as

$$(30) \qquad \Delta(x, x') := \sup_{t \in [x', x]} |Y(t)|.$$

THEOREM 4.10. *For all $x \geq x'$, for all $\varepsilon > 0$ and all integers $n$,*

$$\mathbb{P}(\Delta(x, x') \geq |x - x'|^{1/2+\varepsilon}) = O(|x - x'|^{-n}).$$

This theorem is a consequence of a result on the convergence rates

LEMMA 4.11. *Let $(U_n), n \in \mathbb{N}$, be an i.i.d. sequence of random variables in $\mathbb{R}^2$ and let $(X_i, Y_i)$ be the coordinates of $U_i$. Suppose that $\mathbb{E}X_1 > 0$, $\mathbb{E}Y_1 = 0$ and $\mathbb{E}|U_1|^n < \infty$ for all $n$.*

*Let $t > 0$ and $T_n = \sum_{k=1}^n U_k$. We define*

$$K_t = \inf\left\{ n : \sum_{i=1}^n X_i \geq t \right\} \quad and \quad M_t = \sup_{1 \leq k \leq K_t} \left| \sum_{i=1}^k Y_i \right|.$$

*Then for all $\varepsilon > 0$ and all integers $n$,*

$$\mathbb{P}(M_t \geq t^{1/2+\varepsilon}) = O(t^{-n}).$$

PROOF. Let $c < \mathbb{E}X_1$, $\mathbb{P}(K_t \geq \frac{t}{c}) = P(\sum_{i=1}^{\lfloor t/c \rfloor} X_i \geq t) = O(t^{-n})$, from Theorem 3.1 of [16]. Similarly, from the same theorem,

$$\mathbb{P}\left( \sup_{1 \leq k \leq \lfloor t/c \rfloor} \left| \sum_{i=1}^k Y_i \right| \geq t^{1/2+\varepsilon} \right) = O(t^{-n}).$$

The lemma follows from

$$\mathbb{P}(M_t \geq t^{1/2+\varepsilon}) \leq \mathbb{P}\left( K_t \geq \frac{t}{c} \right) + \mathbb{P}\left( \sup_{1 \leq k \leq \lfloor t/c \rfloor} \left| \sum_{i=1}^k Y_i \right| \geq t^{1/2+\varepsilon} \right). \qquad \square$$

PROOF OF THEOREM 4.10. The $n$th ancestor of $T_0 = (x, 0)$ is $T_n$; let $S_n = T_n - x = -\sum_{k=1}^n U_k$ and if $(X_i, Y_i)$ denotes the coordinates of $U_i$, let

$$K_{x-x'} = \inf\left\{ n : \sum_{i=1}^n X_i \geq x - x' \right\} \quad and \quad M_{x-x'} = \sup_{1 \leq k \leq K_{x-x'}} \left| \sum_{i=1}^k Y_i \right|.$$

We have $\Delta(x, x') = M_{x-x'}$. Since the variables $(U_n), n \in \mathbb{N}$, are not independent, we cannot apply directly Lemma 4.11. Let $(\theta_k), k \in \mathbb{N}$, be the sequence given in Theorem 4.6 and $\tilde{U}_k = \sum_{i=\theta_k+1}^{\theta_{k+1}} U_i$. Due to the regenerative property of the sequence $(\theta_k), k \in \mathbb{N}$, $(\tilde{U}_k), k \in \mathbb{N}^*$ is an i.i.d. sequence; note also that $\mathbb{E}|\tilde{U}_1|^n < \infty$.

Let $\alpha_n = \sup\{k : \theta_{k+1} \le n\}$; we have

$$\left| S_n - \sum_{k=1}^{\alpha_n} \tilde{U}_k \right| \le \sum_{i=1}^{\theta_1} |U_i| + \sum_{i=\theta_{\alpha_n}}^{n} |U_i|.$$

Thus, from Lemma 4.11 applied to $(\tilde{U}_k), k \in \mathbb{N}^*$, we deduce the result for $(U_k), k \in \mathbb{N}$. $\quad \square$

**5. The radial path.** For $x > 0$, let $X = x e_x$ and $R_O(x)$ or $R_O(X)$ denote the path from $X_0 = X = (x, 0)$ to the origin in the RST. This path (a sample of which is depicted in Figure 10) may be seen as a piecewise linear curve in $\mathbb{R}^2$. We denote by $H(X)$ the generation of $T_0$ in the RST. As for the DSF, we will denote by $X_0, X_1, \ldots, X_{H(X)} = O$ the sequence of the successive ancestors of $X_0$ in the RST. Note that the points the sequence $(X_k \cdot e_x), 0 \le k \le H(X)$, is not necessarily a decreasing sequence.

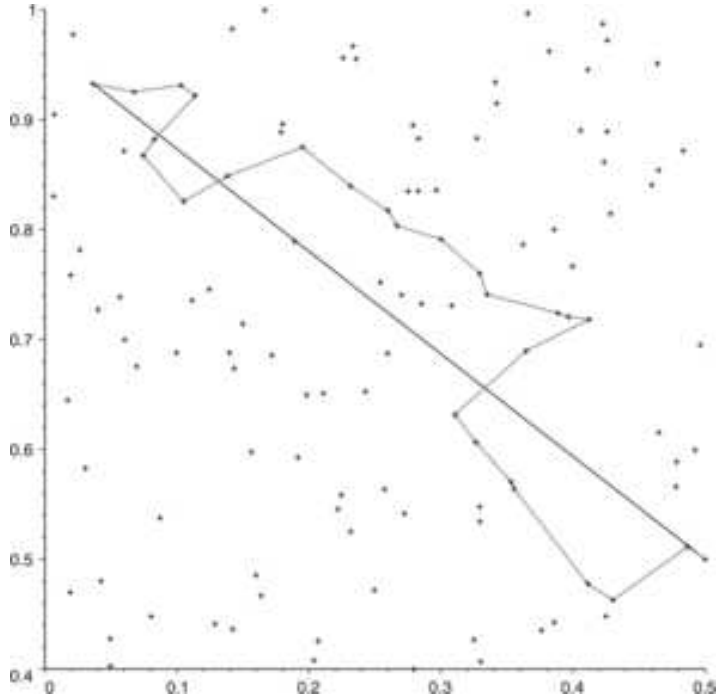

FIG. 10.   *An example of radial path [with its origin in* $(1/2, 1/2)$*]. The initial point is of generation* 26.

5.1. *Maximal deviation of the radial path.* The aim of this section is to bound the tail of the distribution of the maximal deviation of the path

$$\Delta(X) := \sup_{0 \leq k \leq H(X)} |X_k \cdot e_y|. \tag{31}$$

For all $t \in \mathbb{R}$, we define:

$$y(t) := \max(0, \sup\{y : (t, y) \in R_O(X)\}),$$

$$z(t) := \max(0, -\inf\{y : (t, y) \in R_O(X)\}),$$

where the supremum over an empty set is taken to be $-\infty$. Notice that $y(x) = 0$, $y(t) = 0$ for $t \notin [-x, x]$, $\sup_{t \in \mathbb{R}} y(t) = \sup\{Y \cdot e_y : Y \in R_O(X)\}$ and

$$\Delta(X) = \max\left(\sup_{t \in \mathbb{R}} y(t), \sup_{t \in \mathbb{R}} z(t)\right).$$

Now let $\hat{R}(x)$ denote the path from $(x, 0)$ in the directed spanning forest $\hat{\mathcal{T}}_{-e_x}$ with direction $-e_x$, when building this forest on the same point process as the RST, but with all points below the line $\overline{OX}$ removed. The path $\hat{R}(x)$ may be parameterized by the $x$ coordinate of the process

$$\hat{R}(x) = \{(t, \hat{y}(t)), t \leq x\}.$$

The key observation of this subsection is the following:

THEOREM 5.1. *For all $t \leq x$, $y(t) \leq \hat{y}(t)$.*

PROOF. The proof of this result is in two steps.

Let $\mathcal{H}^-(t) = \{Y \in \mathbb{R}^2 : Y \cdot e_x < t\}$ denote the left half-plane $(x < t)$. We first prove that if the two curves share a common point $T$ of the point process, namely, if $(t, y(t)) = (t, \hat{y}(t)) = T$ for some $T$, and if $S \in \mathcal{H}^-(t)$, then $S \cdot e_y \leq S \cdot e_y$, where $S$ (resp. $\hat{S}$) is the ancestor of $T$ in $\mathcal{T}$ (resp. $\hat{\mathcal{T}}_{-e_x}$).

Let $y = y(t) = \hat{y}(t) \geq 0$. If $S = \hat{S}$, then the property is proved.

If $S \neq \hat{S}$ and $S \in L_-$ (where $L_-$ denotes the half-plane $y \leq 0$), the property also holds since $\hat{S} \in L_-^c$ (we recall that the points below $\overline{OX}$ have been removed by construction for $\hat{\mathcal{T}}_{-e_x}$).

If $S \neq \hat{S}$ and $S \in L_-^c$, then $\hat{S} \in \mathcal{H}^-(t) \setminus B(O, |T|)$, and thus, $\hat{S} \in L_-^c \cap m(\mathcal{H}^-(t) \setminus B(O, |T|)m)$, and the property follows in this case too since $S \in B(O, |T|) \cap \mathcal{H}^-(t)$.

To complete the proof, we now show that $R_O(x) \cap \hat{R}(x) \subset N \cup \{O\} \cup \{(x, 0)\}$. Suppose that the two curves intersect in the segment between the vertices $\hat{T}_1 = (\hat{x}_1, \hat{y}_1)$ and $\hat{T}'_1 = (\hat{x}'_1, \hat{y}'_1)$ in $\hat{\mathcal{T}}_{-e_x}$ and the segment between the vertices $T_1 = (x_1, y_1)$ and $T_2 = (x_2, y_2)$ in $\mathcal{T}$. We may assume $\hat{x}_1 \geq \hat{x}_2$ and $x_1 \geq x_2$. Moreover, if $L$ denotes the closed half-plane below the line $(T_1, T_2)$, then $\hat{T}_1 \in L^c$ and $\hat{T}_2 \in L$.

Suppose that $\hat{x}_1 < x_1$. Note that the ancestors of $\hat{T}_1$ (i.e., $\hat{T}_2$) and $T_1$ (i.e., $T_2$) are both in $\mathcal{H}^-(\hat{x}_1)$. From the construction of $\hat{\mathcal{T}}_{-e_x}$, for any point $Y$ of the point process lying in $\mathcal{H}(\hat{x}_1)$, $\hat{T}_2 \in \overline{B(\hat{T}_1, |\hat{T}_1 - Y|)}$, where $\bar{\cdot}$ denotes the closure of a set. Using the fact that $\hat{x}_1 < x_1$, we deduce that $B(\hat{T}_1, |\hat{T}_1 - T_2|) \cap \mathcal{H}^-(\hat{x}_1) \cap L \subset B(T_1, |T_1 - T_2|) \cap \mathcal{H}^-(\hat{x}_1) \cap L$. Now, by construction of the RST, $\overline{B(T_1, |T_1 - T_2|)} \cap \mathcal{H}^-(\hat{x}_1) \cap L \cap N = \{T_2\}$. It follows that $\hat{T}_2 = T_2$.

If $\hat{x}_1 \geq x_1$, then simple geometric arguments show that, necessarily, $\hat{T}_2 = T_1$.  □

COROLLARY 5.2.  *For all $\varepsilon > 0$ and all integers $n$,*

$$\mathbb{P}(\Delta(X) \geq |X|^{1/2+\varepsilon}) = O(|X|^{-n}).$$

PROOF.  Let $x > 0$, $X = (x, 0)$ and $\tilde{R}(x)$ be the path from $X$ in the directed spanning forest with direction $-e_x$, obtained when using the same point process as in the RST but when removing all points above the line $\overline{OX}$. The path $\tilde{R}(x)$ may be parameterized as $\tilde{R}(x) = \{(t, -\tilde{z}(t)), t \leq x\}$. By symmetry, the processes $\tilde{z}$ and $\hat{y}$ have the same distribution and from Theorem 5.1, $z(t) \leq \tilde{z}(t)$. Therefore, it is sufficient to prove the following:

$$\mathbb{P}\left(\sup_{t \in [-x, x]} \hat{y}(t) \geq x^{1/2+\varepsilon}\right) = O(x^{-n}).$$

A similar statement for the path $R(x)$ in the directed chain built on the whole Poisson point process was proved in Theorem 4.10. The difference between $R(x)$ and $\hat{R}(x)$ appears on the vicinity of the line $\overline{OX}$. Thus, the supremum should not be much affected.

To formalize this idea, consider the event, for $x > 0$, $\gamma > 0$,

$$A_{\gamma,x} = \{\exists Y \in \mathbb{R}^2, \text{ with } |Y| \leq 2x \text{ and } |Y - \mathcal{A}(Y)| \geq x^\gamma\}.$$

Following the proof of Lemma 5.2 of [17], we have

$$(32) \qquad \mathbb{P}(A_{\gamma,x}) \leq C_1 \exp(-C_0 x^{2\gamma}).$$

Fix $0 < \varepsilon'' < \varepsilon' < \varepsilon$ and $\gamma \leq \frac{1}{2}$. On $A_{\gamma,x}^c$, if $\sup_{t \in [-x,x]} \hat{y}(t) \geq x^{1/2+\varepsilon}$, then for $x$ large enough, there exists at least one point of $\hat{R}(x) \cap N$ with its $y$ coordinate in the interval $[x^{1/2+\varepsilon'}, 2x^{1/2+\varepsilon'}]$. Let $X_0 = (t_0, y_0)$ be the rightmost such point. Note that $|X_0| \geq x^{1/2+\varepsilon'} \geq \sqrt{x}$. The path $\hat{R}(x)$ coincides with the path $R(X_0)$ on the interval $[-x, X_0]$, provided the maximal deviation of $R(X_0)$ is less than $x^{1/2+\varepsilon''}$; indeed, the infimum of $R(X_0)$ will then be lower bounded by $y_0 - x^{1/2+\varepsilon''} > 0$. The event $\{R(X_0) \geq x^{1/2+\varepsilon''}\}$ is a subset of $\tilde{A}_{\varepsilon'',x} = \bigcup_{Y \in N \cap B(O,x) \setminus B(O,\sqrt{x})} \{\Delta(|Y|, -|Y|)(R(Y)) \geq x^{1/2+\varepsilon''}\}$.

From Theorem 4.10 there exits a constant $C$ (depending on $\varepsilon''$) such that $\mathbb{P}(\Delta(|Y|, -|Y|)(R(Y)) \geq |Y|^{1/2+\varepsilon''}) \leq C|Y|^{-n}$. Now, letting $l > e$, we get

$$\mathbb{P}(\tilde{A}_{\varepsilon'',x}) \leq \mathbb{P}\left( \bigcup_{Y \in N \cap B(O,x) \setminus B(O,\sqrt{x})} \{\Delta(|Y|, -|Y|)(R) \geq |Y|^{1/2+\varepsilon''}\} \right)$$

$$\leq \mathbb{P}(N(B(O,x)) \geq l\pi x^2) + l\pi x^2 C|x|^{-n/2}.$$

Thus, $\mathbb{P}(\tilde{A}_{\varepsilon'',x}) = O(|x|^{-n})$ for all integer $n$. Using this last remark, the conclusion follows from equation (32) and the inequality

$$\mathbb{P}\left( \sup_{t \in [-x,x]} \hat{y}(t) \geq x^{1/2+\varepsilon} \right) \leq \mathbb{P}(A_{\gamma,x}) + \mathbb{P}(\tilde{A}_{\varepsilon'',x}). \qquad \square$$

For $X \in N$, let $R_{\text{out}}(X)$ be the set of offspring of $X$ in the RST, namely, the set $X' \in N$ such that $X \in R_O(X')$. We now state a definition introduced in [17].

DEFINITION 5.3.  Let $f$ be a positive function on $\mathbb{R}_+$. A tree is said to be $f$-straight at the origin, if for all but finitely many vertices

$$R_{\text{out}}(X) \subset C(X, f(X)),$$

where for all $X \in \mathbb{R}^2$ and $\varepsilon \in \mathbb{R}_+$, $C(X, \varepsilon) = \{Y \in \mathbb{R}^2 : \theta(X, Y) \leq \varepsilon\}$ and $\theta(X, Y)$ is the absolute value of the angle (in $[0, \pi]$) between $X$ and $Y$.

THEOREM 5.4.  *The number of points of $N$ such that $\Delta(X) \geq |X|^{1/2+\varepsilon}$ is a.s. finite. The RST is $f$-straight at $O$ for $f(X) = |X|^{-1/2+\varepsilon}$ for all $\varepsilon > 0$.*

PROOF.  We first prove that the number $K$ of points $T_n$ of $N$ such that $\Delta(T_n) \geq |T_n|^{1/2+\varepsilon}$ is finite.

From Corollary 5.2, for all $X \in \mathbb{R}^2$, $\mathbb{P}(\Delta(X) \geq |X|^{1/2+\varepsilon}) \leq 1 \wedge C|X|^{-3}$, where $C$ depends on $\varepsilon$. From Campbell's formula,

$$\mathbb{E}K = \mathbb{E} \sum_{T_n \in N} \mathbb{1}(\Delta(T_n) \geq |T_n|^{1/2+\varepsilon})$$

$$= 2\pi \int_0^\infty \mathbb{P}(\Delta((x,0)) \geq x^{1/2+\varepsilon}) x \, dx$$

$$\leq 2\pi \int_0^\infty 1 \wedge C x^{-2} \, dx < \infty.$$

The rest of the proof uses the same argument as Theorem 2.6, Lemma 2.7 of [17] (with 3/4 replaced with 1/2).  $\square$

$f$-straight trees have a simple topology described by Proposition 2.8 of [17] and restated in Theorem 2.1.

5.2. *Law of large numbers on the radial path.* The aim of this section is to prove Theorem 2.5. We use the notation introduced at the beginning of Section 5. Let $X = (x, 0)$ and $X_0 = X, X_1, \ldots, X_{H(X)} = O$ be the successive ancestors of $X$ in the RST.

The edge between $X_i$ and $X_{i-1}$ is $U_i = X_{i-1} - X_i$. Letting $U_0 = -X_0$, we have

$$X_n = -\sum_{l=0}^{n} U_l \quad \text{and} \quad X_n - X_k = -\sum_{l=k+1}^{n} U_l.$$

We define $D_0 := B(X_0, |X_0|)$ and for $n \geq 1$,

$$(33) \qquad D_n := B(X_n, |X_n|) \cap \left\{ \bigcup_{k=0}^{n-1} B(X_n - X_k, |U_{k+1}|) \right\}^c.$$

The density of $U_{n+1}$ given $(U_0, \ldots, U_n)$ is

$$\mathbb{1}(r < |X_n|) \frac{d}{dr} \nu_2(B(O, r) \cap D_n) e^{-\nu_2(B(O, r) \cap D_n)} \frac{\mathbb{1}_{(r, \theta) \in D_n} r \, dr \, d\theta}{\nu_1(C(O, r) \cap D_n)} + c \delta_{-X_n}(r, \theta),$$

where $c$ is a normalizing constant.

Letting $g$ be a measurable $\mathbb{R}^2 \mapsto \mathbb{R}$ function, we suppose that, for all $Y \in \mathbb{R}^2$, $|g(Y)| \leq \max(C, |Y|^\alpha)$ some positive constants $C$ and $\alpha$. We define the statistical average of $g$ on $R_O(X)$ as

$$S_{H(X)} = \sum_{n=1}^{H(X)} g([X_{n-1} - X_n]_{X_{n-1}}),$$

where $[U]_V$ is the vector $U$ rotated by an angle $-\theta$ and $V = r \cos \theta \cdot e_x + r \sin \theta \cdot e_y$ ($[U]_V$ is the vector $U$ expressed in the local coordinates of $V$).

We will first prove that

$$(34) \qquad \lim_{|X| \to +\infty} \frac{S_{H(X)}}{H(X)} = \pi(g),$$

where $\pi$ is defined by equation (26).

Results on the matter will follow from the corresponding results on the directed path. The proof of the theorem is in three steps.

5.2.1. *Step I: pseudo-regenerative times.* In this part we associate to the RST a sequence analog to the sequence given by Theorem 4.6 for the DSF.

Let $\mathcal{F}_n$ be the $\sigma$-field generated by $X_0, U_1, \ldots, U_n$. Clearly, $X_n \in \mathcal{F}_n$. We build a coupling on the process $\{U_n\}$ in the same vein as in [3] for Markov chains. As in the DSF, we define

$$\tau_{n+1} = \inf \left\{ m > \tau_n : B(O, |X_m|) \cap \left\{ \bigcup_{k=\tau_n - 1}^{m-2} B(X_m - X_k, |U_{k+1}|) = \varnothing \right\} \right\},$$

with $\tau_0 = 1$ and $\Phi_n = \{0, X_{\tau_{n-1}} - X_{\tau_{n-1}+1}, \ldots, X_{\tau_{n-1}} - X_{\tau_n}\}$. For $n$ large enough, $\Phi_n = \{0\}$. Let $\widetilde{\Phi}_n$ denote the $\Phi_n$ vector in the following basis: the origin is in $X_{\tau_{n-1}}$ and the basis vectors are $(e_{X_{\tau_{n-1}}}, e_{\hat{X}_{\tau_{n-1}}})$, where $\hat{X}_{\tau_{n-1}}$ is the point of coordinates $(-X^{(2)}_{\tau_{n-1}}, X^{(1)}_{\tau_{n-1}})$.

We recall some notation introduced in the proof of Theorem 4.4. For $\widetilde{\Phi} = \{0, Y_1, \ldots, Y_{\tau(\Phi)}\}$ in $\mathcal{N}_{\mathcal{X}}$ and $\tau(\Phi) \geq 1$, we define $(R(\widetilde{\Phi}), \Theta(\widetilde{\Phi}))$ as the coordinates of $Y_{\tau(\Phi)}$. The set $\mathcal{S} \in \mathcal{B}(\mathcal{N}_{\mathcal{X}})$ is defined by

$$(35) \qquad \mathcal{S} = \{0\} \cup \{\Phi \in \mathcal{N}_{\mathcal{X}} \setminus \{0\} : |\Theta(\Phi)| \leq \pi/6, R(\Phi) \leq 1\}.$$

Let $\mathcal{G}$ denote the sub $\sigma$-algebra of $\mathcal{B}(\mathcal{N}_{\mathcal{X}})$ generated by the sets $\{\tau(\Phi) = 1\} \cap \{X_1(\Phi) \in B\}, B \in C$, where $C = \{(r, \theta) \in \mathbb{R}^2 : \theta \in [-\pi/3, \pi/3], r \leq 1\}$. We have the following lemma:

LEMMA 5.5.    *There exist positive constants $\lambda$ and $g$, and a measure $\nu$ with support on $\mathcal{G}$, such that, for any $x_0$ large enough and for all $n$:*

(i) $\mathbb{P}(\widetilde{\Phi}_{n+1} \in \mathcal{S} | \mathcal{F}_{\tau_n}) \geq \lambda > 0$;

(ii) *if $|X_{\tau_n}| \geq x_0$ and $\widetilde{\Phi}_n \in \mathcal{S}$, then $\mathbb{P}(\widetilde{\Phi}_{n+1} \in A | \mathcal{F}_{\tau_n}) \geq g\nu(A)$.*

PROOF.    Define $c_\alpha(R) := \{X = (r, \theta) \in \mathbb{R}^2 : \theta \in (\alpha, 0], r \leq R\}$ for $\alpha < 0$, and $c_\alpha(R) := \{X = (r, \theta) \in \mathbb{R}^2 : \theta \in [0, \alpha), r \leq R\}$ for $\alpha > 0$. We may prove as in Lemma 4.2 that the RST is such that either $c_{\pi/6}(|X_n|) \subset \widetilde{D}_n$ or $c_{-\pi/6}(|X_n|) \subset \widetilde{D}_n$, where $\widetilde{D}_n$ is the set defined by equation (33) re-expressed in the above referential. Hence,

$$\mathbb{P}(\widetilde{U}_{n+1} \in c_{\pi/6}(1) \cup c_{-\pi/6}(1) | \mathcal{F}_n) \geq \frac{1}{6}(1 - e^{-\pi/12}) = \lambda > 0.$$

The proof of property (i) is then similar to that of (23) in Lemma 4.5 for the RST case.

For $x_0$ large enough, if $|X_{\tau_n}| \geq x_0$ and $\widetilde{\Phi}_n \in \mathcal{S}$, the cone $C$ is included in $\widetilde{D}_{\tau_n}$. As in the proof of Theorem 4.4, we can check that $\nu = P_0^1|_{\mathcal{G}}$ satisfies property (ii) for a suitable choice of $g$, where $P_0^1$ was defined in Theorem 4.4. $\square$

Fix $x_0$ and define

$$(36) \qquad\qquad n(X) := \sup\{n : |X_n| \geq x_0\}.$$

Following the proof of the regeneration lemma [3], we can show the existence of a sequence of stopping times $N_k$ (on an enlarged probability space) such that $N_0 = 0$ and such that, for all $A \in \mathcal{B}(\mathcal{N}_{\mathcal{X}})$,

$$\mathbb{P}(\Phi_{n+1} \in A, N_k = n, \tau_n \leq n(X)) = \nu(A)\mathbb{P}(N_k = n, \tau_n \leq n(X)).$$

From property (i) and (ii), $\{N_k\}$ is a renewal process and it can be built such that $N_{k+1} - N_k$ is a geometric random variable of parameter $\lambda g > 0$.

Let

$$K'(X) := \sup\{k : |X_{\tau_{N_k}}| \geq x_0\}.$$

For $1 \leq k \leq K'(X)+1$, $(\widetilde{\Phi}_{N_k+1}, \ldots, \widetilde{\Phi}_{N_{k+1}})$ depends on the past only through $|X_{\tau_{N_{k-1}}}|$. Thus, we have found some homogeneous regenerative structure on the nonhomogeneous process $(U_n)$.

Let

(37) $$\theta_k = \tau_{N_{k-1}}$$

and

(38) $$K(X) = \tau_{N_{K'(X)}}.$$

Note that Lemma 4.3 is still valid for the RST. Indeed, this lemma relies on the bounds $\tilde{L}_i$ and $\tilde{P}_i$ which also hold for the RST. Finally, we have the following decomposition:

LEMMA 5.6. *For all $x_0$ large enough, there exists an increasing sequence of finite stopping times $\{\theta_k\}, k \in \mathbb{N}$, (on an enlarged probability space) such that $\theta_0 = 0$, for all $k \geq 0$, $\mathbb{P}(\theta_{k+1} - \theta_k > n | \mathcal{F}_{\theta_k}) \leq C_1 \exp(-C_0 n)$, for some positive constants and for $0 \leq k \leq K(X)$, $(U_{\theta_k+1}, \ldots, U_{\theta_{k+1}})$ depends on $\mathcal{F}_{\theta_k}$ only through $|X_{\theta_k}|$.*

LEMMA 5.7. *Almost surely, $H(X)$ and $K(X)$ tend toward infinity with $|X|$.*

PROOF. For $\beta > 0$, $x > 0$, let $A_{\beta,x} = \{\exists Y \in N,$ with $|Y| \leq 2x$ and $|Y - \mathcal{A}(Y)| \geq x^\beta\}$. The proof of Lemma 5.2 of [17] implies that $\mathbb{P}(A_{\beta,x}) = O(x^{-q})$ for all $q$.

On $A_{\beta,|X|}^c$, $H(X) \geq \lfloor |X|^{1-\beta} \rfloor$. Let $M$ be the number of points in $N$ such that $H(X) \leq |X|^{1-\beta}$. From Campbell's formula, $\mathbb{E}M \leq 2\pi \int_0^\infty \mathbb{P}(A_{\beta,x}) x \, dx < \infty$. Hence, $M$ is finite almost surely.

To prove that $K(X)$ is finite, notice that

$$\mathbb{P}(\theta_{k+1} - \theta_k > n | \mathcal{F}_{\theta_k}) \leq C_1 \exp(-C_0 n),$$

which implies that a.s. $\limsup_k \theta_k / k < C_1 < \infty$. Hence, a.s. for $k$ large enough, $|X_{\lfloor C_1 k \rfloor}| \geq x_0$ implies that $K(X) \geq k$. Pick $k = \lfloor |X|^{1-\beta} \rfloor$ and the proof is similar to the proof for $H(X)$. $\quad \square$

5.2.2. *Step II: identification of the limit.* For $0 \le k \le K(X)$, let $V_n = [U_n]_{X_{n-1}}$ and

$$\tilde{S}_k(X) = \sum_{n=\theta_k+1}^{\theta_{k+1}} g(V_n).$$

The next step is to identify $\lim_{|X|\to\infty} \mathbb{E}\tilde{S}_k(X)$. Notice that in the DSF with direction $-e_x$, the same construction can be done on the path starting from $O$; we will denote similarly $U_n^{DSF}$, $\theta_k^{DSF}$ and $\tilde{S}_k^{DSF} = \sum_{n=\theta_k^{DSF}+1}^{\theta_{k+1}^{DSF}} g(U_n^{DSF})$. These sequences are i.i.d. sequence and, since $X = xe_x$, Theorem 4.7 gives $\mathbb{E}\tilde{S}_k^{DSF} = \pi(g)\overline{\theta}$, where $\overline{\theta} = \mathbb{E}(\theta_{k+1}^{DSF} - \theta_k^{DSF})$.

From Lemma 5.6, using the Cauchy–Schwarz inequality, we get

$$\begin{aligned}
(39) \qquad \mathbb{E}\tilde{S}_1(X) &= \sum_{k=0}^{\infty} \mathbb{E}\left( \mathbb{1}(\theta_1 = k) \sum_{n=0}^{k} g(V_n(X)) \right) \\
&\le \sum_{n=0}^{\infty} C_1 \exp(-C_0 n) \sqrt{\mathbb{E}g^2(V_n(X))}.
\end{aligned}$$

Since $g(V_n(X))$ is stochastically dominated by $\max(C, Z^\alpha)$, where $\mathbb{P}(Z > t) \le e^{-\pi t^2/12}$, it follows from Theorem 2.4 that

$$\lim_{|X|\to\infty} \mathbb{E}g(V_n(X)) = \mathbb{E}g(U_n^{DSF}).$$

Since $\mathbb{P}(|U_{n+1}(X)| > t|\mathcal{F}_n) \le e^{-\pi t^2/12}$, $\mathbb{E}g^2(V_n(X))$ is uniformly bounded. So, from the dominated convergence theorem and equation (39), we deduce that $\lim_{|X|\to+\infty} \mathbb{E}\tilde{S}_1(X) = \pi(g)\overline{\theta}$.

Similarly, the law $\tilde{S}_k(X)$ depends on $X$ and $k$ only through $|X_{\theta_k}|$, which tends toward infinity as $|X|$ tends toward infinity; it follows a.s.

$$(40) \qquad \lim_{|X|\to+\infty} \mathbb{E}(\tilde{S}_k(X)|\mathcal{F}_{\theta_k}) = \pi(g)\overline{\theta}.$$

5.2.3. *Step III: convergence.* We may write

$$\frac{S_{H(X)}}{H(X)} = \frac{1}{H(X)} \sum_{k=1}^{K(X)} \tilde{S}_k(X) + \frac{1}{H(X)} \sum_{\theta_{K(X)+1}+1}^{H(X)} g(V_n).$$

There is almost surely a finite number of edges in $B(O, x_0)$; thus, from Lemma 5.7, the second term tends a.s. to 0. From equation (40), we also obtain that a.s. $\lim_{|X|\to\infty} 1/K(X) \sum_{k=1}^{K(X)} \mathbb{E}(\tilde{S}_k(X)|\mathcal{F}_{\theta_k}) = \pi(g)\overline{\theta}$.

Moreover, using the same arguments as in equation (39), we get $\mathbb{E}(|\tilde{S}_k(X)|^2) < M < \infty$. We can then apply Theorem VII.9.3 of [11] which gives

$$\lim_{|X|\to\infty} \frac{1}{K(X)} \sum_{k=1}^{K(X)} \tilde{S}_k(X) = \lim_{|X|\to+\infty} \frac{1}{K(X)} \sum_{k=1}^{K(X)} \mathbb{E}(\tilde{S}_k(X)|\mathcal{F}_{\theta_k}) = \pi(g)\overline{\theta}, \qquad \text{a.s.}$$

For $g = 1$, we deduce that a.s. $H(X)/K(X)$ tends to $\overline{\theta}$. This ends the proof of equation (34).

To complete the proof of Theorem 2.5, it remains to show that if $g$ is continuous, $\lim_{|X|\to+\infty} 1/H(X) \sum_{k=1}^{H(X)} g(X_{k-1} - X_k) = \pi(g)$. This is a consequence of Theorem 5.4 and equation (34). Indeed, a.s. for all $k$, as $|X|$ tends toward infinity, $[X_{k-1} - X_k]_{X_{k-1}}$ tends toward $X_{k-1} - X_k$. If $g$ is continuous, we deduce that $1/H(X) \sum_{k=1}^{H(X)} g(X_{k-1} - X_k)$ and $1/H(X) \sum_{k=1}^{H(X)} g([X_{k-1} - X_k]_{X_{k-1}})$ have the same limit.

COROLLARY 5.8.    *The following a.s. limits hold*:

$$\lim_{|X|\to+\infty} \frac{H(X)}{|X|} = \frac{1}{p},$$

$$\lim_{|X|\to+\infty} \frac{1}{|X|} \sum_{k=0}^{H(X)-1} |T_{k+1} - T_k|^\alpha = \frac{l_\alpha}{p},$$

*where $p$ and $l_\alpha$ are defined in Corollary 4.8.*

5.3. *Shape theorem.* We define $G_k := |\mathcal{T}(k)|$; $G_k$ is the size of the ball of center $O$, and radius $k$ for the graph-distance on the RST.

The main aim of this section is to prove Theorem 2.2, and in particular the fact that $G_k/k^2$ a.s. tends to a constant when $k$ tends to $\infty$. In the literature this constant is known as the *volume growth*.

The intuition behind Theorem 2.2 is as follows: from the results of Section 5.2, a point $k$ hops away from the origin is asymptotically at Euclidean distance $d_k \sim kp$ from the origin. The disc of radius $d_k$ contains $\pi d_k^2$ vertices asymptotically.

In order to prove Theorem 2.2, we need an estimate of the tail of $H(X)$ around its mean. In Figure 11 is shown a realization of $\mathcal{T}(100)$.

The proof of the next theorem is the heart of the proof of Theorem 2.2.

THEOREM 5.9.    *For all $q < p$, there exists positive constants $C_0$ and $C_1$:*

$$\mathbb{P}\left(H(X) > \frac{|X|}{q}\right) \le C_1 \exp(-C_0|X|).$$

*Similarly, for $q > p$,*

$$\mathbb{P}\left(H(X) < \frac{|X|}{q}\right) \le C_1 \exp(-C_0|X|).$$

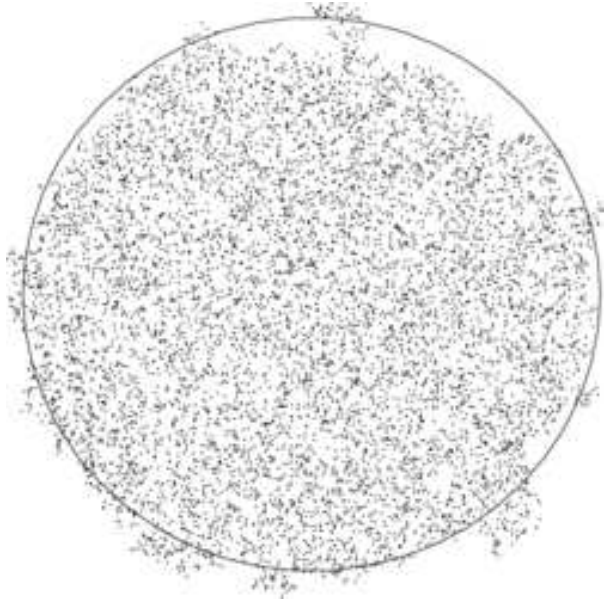

Fig. 11. *The points of $\mathcal{T}(100)$ and the disc of radius $100p$.*

The proof is again in three steps.

Let $X_0 = X, X_1, \ldots, X_{H(x)} = O$, the successive ancestors of $X$ in the radial spanning tree. For $1 \le n \le H(X)$, we define the progress

$$P_n(X) = |X_{n-1}| - |X_n|$$

and for $n > H(X)$, $P_n(X) = 0$. By definition, for all integers $m$,

$$\mathbb{P}(H(X) > m) = \mathbb{P}\left(\sum_{n=1}^{m} P_n(x) < |X|\right).$$

In order to derive the tail bounds of the theorem, we will analyze the asymptotic behavior of $P_n(X)$.

5.3.1. *Step I: pseudo regenerative times.* We will use the notation of Section 5.2 and the sequence of stopping times $0 = \theta_0, \theta_1, \ldots$ given by Lemma 5.6, which are such that, for $0 \le k \le K(X)$, $(P_{\theta_k+1}(X), \ldots, P_{\theta_{k+1}}(X))$ depends of $\mathcal{F}_{\theta_k}$ only through $|X_{\theta_k}|$.

Let

$$\widehat{P}_k(X) = \sum_{n=\theta_k+1}^{\theta_{k+1}} P_n(X).$$

From equation (40),

$$(41) \qquad \lim_{|X| \to \infty} \mathbb{E}\,\widehat{P}_k(X) = p\overline{\theta} = p/\delta,$$

where $\overline{\theta} = \delta^{-1} = \lim_{|X| \to \infty} \mathbb{E}\theta_1$.

LEMMA 5.10.   *There exists positive constants $C_0$, $C_1$ such that, for all $k$,*
$$\mathbb{P}(\widehat{P}_k(X) \geq t | \mathcal{F}_{\theta_k}) \leq C_1 \exp(-C_0 t).$$

PROOF.   We will prove this result for $k = 1$; the proof extends to all $k$ with minor changes. From Lemma 5.6, the sequence $(\theta_k)$ is such that, for some constants $C_0$, $C_1$, $\mathbb{P}(\theta_{k+1} - \theta_k \geq K | \mathcal{F}_{\theta_k}) \leq C_1 \exp(-C_0 K)$ for all $K > 0$. Let $s > 0$. From the Cauchy–Schwarz inequality,

$$\begin{aligned}
\mathbb{E}e^{s\widehat{P}_1(X)} &= \sum_{K=1}^{\infty} \mathbb{E}\mathbb{1}(\theta_1 = K)e^{s\sum_{n=1}^{K} P_n(X)} \\
&\leq \sum_{K=1}^{\infty} \sqrt{\mathbb{P}(\theta_1 = K)\mathbb{E}e^{2s\sum_{n=1}^{K} P_n(X)}} \\
&\leq \sum_{K=1}^{\infty} C_1 \exp(-C_0 K)\sqrt{\mathbb{E}e^{2s\sum_{n=1}^{K} P_n(X)}}.
\end{aligned}$$
(42)

We have already noticed that

$$\mathbb{P}(P_{n+1}(X) \geq t | \mathcal{F}_n) \leq \exp\left(-\frac{\pi}{12}t^2\right),$$

$$\mathbb{E}(P_{n+1}(X) | \mathcal{F}_n) \leq \mu.$$

Calculations based on Chernoff's inequality lead to

$$\mathbb{E}e^{2s\sum_{n=1}^{K} P_n(X)} \leq e^{sK\mu}(1 + C_1 s)^K e^{KC_1^2 s^2} \leq C_1 e^{C_1' sK},$$

for all $s \leq 1$. For $s$ sufficiently small, $sC_1' < C_0$, where $C_0$ is the constant in (42). For such a choice of $s$, we get $\mathbb{E}e^{s\widehat{P}_1(X)} \leq C_1$ from equation (42). The lemma then follows from Chernoff's inequality.   □

5.3.2. *Step II: case $q < p$.*   Let $q < p' < p$ and $\varepsilon > 0$ be such that $q < \frac{\delta}{\delta + \varepsilon}p' < p$, with $\delta = \overline{\theta}^{-1}$. From equation (41), there exists $x_0$ such that for $|X| > x_0$, $\mathbb{E}(\widehat{P}_k(X)||X_{\theta_k}| \geq x_0) > \frac{p'}{\delta + \varepsilon}$:

$$\begin{aligned}
\mathbb{P}(H(X) &> m) \\
&= \mathbb{P}\left(\sum_{n=1}^{m} P_n(X) < |X|\right) \\
&\leq \mathbb{P}\left(\sum_{n=1}^{m-l} P_n(X) < |X| - x_0\right) + \mathbb{P}(N(B(O, x_0) > l)) \\
&\leq \mathbb{P}\left(\sum_{n=1}^{m-l} P_n(X) < |X| - x_0\right) + \exp\left(-l\ln\frac{l}{e\pi x_0^2}\right),
\end{aligned}$$
(43)

where we have used the classical inequality

$$\mathbb{P}(N(B(O, x_0)) > l) \leq \exp\left(-l \ln \frac{l}{e\pi x_0^2}\right)$$

(Lemma 11.1.1 of [28]).

Let $n(x)$ and $K(x)$ be defined as in (36) and (38), respectively. Let $Z_k$ be the conditional random variable $\widehat{P}_k(X)$ given $K(X) \geq k$. We have $\mathbb{E}(Z_k|\mathcal{F}_{\theta_k}) > \frac{p'}{\delta+\varepsilon}$. Note that, for all integers $j$,

$$
\begin{aligned}
&\mathbb{P}\left(\sum_{n=1}^{j} P_n(X) < |X| - x_0\right) \\
&\quad = \mathbb{P}\left(\sum_{n=1}^{j} P_n(X) < |X| - x_0, n(X) > j\right) \\
&\quad \leq \mathbb{P}\left(\sum_{k=1}^{\lfloor(\delta+\varepsilon)j\rfloor} \widehat{P}_k(X) < |X| - x_0\right) + \mathbb{P}(\theta_{\lfloor(\delta+\varepsilon)j\rfloor} < j) \\
&\quad \leq \mathbb{P}\left(\sum_{k=1}^{\lfloor(\delta+\varepsilon)j\rfloor} Z_k < |X|\right) + \mathbb{P}\left(\sum_{k=1}^{\lfloor(\delta+\varepsilon)j\rfloor} \theta_k - \theta_{k-1} < j\right).
\end{aligned}
\tag{44}
$$

Indeed, $\{\sum_{k=1}^{\lfloor(\delta+\varepsilon)j\rfloor} \widehat{P}_k(X) < |X| - x_0\}$ is a subset of $\bigcup_{1 \leq k \leq (\delta+\varepsilon)j}\{K(X) \geq k\}$.

Now, let $x, j$ such that $p'j - x > \alpha j$ for some $\alpha > 0$:

$$
\begin{aligned}
&\mathbb{P}\left(\sum_{k=1}^{\lfloor(\delta+\varepsilon)j\rfloor} Z_k < |X|\right) \\
&\quad = P\left(\frac{1}{(\delta+\varepsilon)j} \sum_{k=1}^{\lfloor(\delta+\varepsilon)j\rfloor} (Z_k - \mathbb{E}(Z_k|\mathcal{F}_{N_k}))\right. \\
&\qquad\qquad \left. < \frac{|X|}{(\delta+\varepsilon)j} - \frac{1}{(\delta+\varepsilon)j} \sum_{k=1}^{\lfloor(\delta+\varepsilon)j\rfloor} \mathbb{E}(Z_k|\mathcal{F}_{N_k})\right) \\
&\quad \leq P\left(\frac{1}{(\delta+\varepsilon)j} \sum_{k=1}^{\lfloor(\delta+\varepsilon)j\rfloor} (Z_k - \mathbb{E}(Z_k|\mathcal{F}_{N_k})) < \frac{|X|}{(\delta+\varepsilon)j} - \frac{p'}{q+\varepsilon}\right) \\
&\quad \leq P\left(\left|\frac{1}{(\delta+\varepsilon)j} \sum_{k=1}^{\lfloor(\delta+\varepsilon)j\rfloor} (Z_k - \mathbb{E}(Z_k|\mathcal{F}_{N_k}))\right| > \frac{jp' - |X|}{(\delta+\varepsilon)j}\right).
\end{aligned}
$$

The process $(Z_k - \mathbb{E}(Z_k|\mathcal{F}_{N_k}))_k$ satisfies the hypothesis of Lemma A.1 (in the Appendix) with $t_0 = \frac{\alpha}{\delta + \varepsilon}$. We deduce that

$$(45) \qquad \mathbb{P}\left( \sum_{k=1}^{\lfloor (\delta + \varepsilon) j \rfloor} Z_k < |X| \right) \leq C_1 e^{-C_0 (jp' - |X|)}.$$

Similarly, $\theta_k - \theta_{k-1}$ satisfies the same light tail hypothesis as $Z_k$ and for $x_0$ large enough, $\mathbb{E}(\theta_k - \theta_{k-1}||X_{\theta_{k-1}}| \geq x_0) > t_0 > 1/(\delta + \varepsilon)$. We get

$$(46) \qquad \mathbb{P}\left( \sum_{k=1}^{\lfloor (\delta + \varepsilon) j \rfloor} \theta_k - \theta_{k-1} < j \right) \leq C_1 e^{C_0 j}.$$

Now let $m = \lfloor \frac{x}{q} \rfloor$, $l = \lfloor \frac{\varepsilon |X|}{2q\delta} \rfloor$ and $j = m - l$. For $|X|$ large enough, the inequality $p'j - |X| > \alpha j$ holds for some positive $\alpha$. Putting together equations (43), (44), (45) and (46), we obtain that, for $|X|$ large enough,

$$\mathbb{P}\left( H(X) > \frac{|X|}{q} \right) \leq C_1 e^{-C_0 |X| \ln |X|} + C_1 e^{C_0 |X|} + C_1 e^{C_0 |X|} = C_1 e^{C_0 |X|}.$$

Since $\mathbb{P}(H(X) > \frac{|X|}{q}) \leq 1$, by increasing arbitrarily $C_1$, we get the above inequality for all $X$.

5.3.3. *Step III: case $q > p$.* The bound $\mathbb{P}(H(X) < \frac{|X|}{q})$ for $q > p$ is obtained with the same type of decomposition:

$$\mathbb{P}(H(X) \leq m)$$
$$= P\left( \sum_{n=1}^{m} P_n(X) = |X| \right)$$
$$\leq P\left( \sum_{n=1}^{m} P_n(X) \geq |X| - x_0 \right)$$
$$\leq P\left( \sum_{k=1}^{\lfloor (\delta - \varepsilon) m \rfloor} \tilde{P}_k(X) \geq |X| - x_0 \right) + \mathbb{P}(\theta_{\lfloor (\delta - \varepsilon) m \rfloor} > m)$$
$$\leq P\left( \sum_{k=1}^{\lfloor (\delta - \varepsilon) m \rfloor} Z_k \geq |X| - x_0 \right) + \mathbb{P}\left( \sum_{k=1}^{\lfloor (\delta - \varepsilon) m \rfloor} \theta_k - \theta_{k-1} > m \right)$$
$$\leq C_1 e^{-C_0 (|X| - x_0 - mp')} + C_1 e^{-C_0 m},$$

for $p' > p$ and all $X$, $m$ such that $|X| - mp' > \alpha m$, for some $\alpha > 0$. Taking $m = \lfloor |X|/q \rfloor$ and $q > p' > p$, we obtain the upper bound for $\mathbb{P}(H(X) < |X|/q)$ and we are done with the proof of Theorem 5.9.

Note that the statement of the theorem is not tight since we do not bound the tail of $|H(X) - |X|/p|$. It seems possible to derive some concentration inequalities.

REMARK 5.11.   Letting $S(X) = \frac{1}{H(X)} \sum_{n=1}^{H(X)} f(U_n)$, the same type of theorem holds for any function $f$ which is growing to infinity slower than exponentially, with $1/p$ replaced by $\pi(f)/p$.

5.3.4. *Proof of Theorem* 2.2.   We begin with the proof of equation (2). Notice that

$$G_k = \sum_{T \in N} \mathbb{1}(H(T) \leq k).$$

Letting $\varepsilon \in (0, 1)$, we may write

$$\frac{|G_k - N(B(O, pk))|}{k^2}$$

$$\leq \frac{1}{k^2} \Bigg( \sum_n \mathbb{1}(X_n \notin B(O, pk) \cap H(X_n) \leq k)$$

$$+ \sum_n \mathbb{1}(X_n \in B(O, pk) \cap H(X_n) > k) \Bigg)$$

$$\leq \frac{1}{k^2} \sum_n \mathbb{1}(X_n \notin B(O, (1+\varepsilon)pk) \cap H(X_n) \leq k)$$

$$+ \frac{N(B(O, (1+\varepsilon)pk) \setminus B(O, (1-\varepsilon)pk))}{k^2}$$

$$+ \frac{1}{k^2} \sum_n \mathbb{1}(X_n \in B(O, (1-\varepsilon)pk) \cap H(X_n) > k)$$

$$\leq I_k + J_k + L_k.$$

From Campbell's formula and using Theorem 5.9,

$$\mathbb{E}(I_k) = \frac{2\pi}{k^2} \int_{(1+\varepsilon)pk}^{\infty} \mathbb{P}(H(xe_x) \leq k)x \, dx$$

$$\leq \frac{2\pi}{k^2} \int_{(1+\varepsilon)pk}^{\infty} C_1 e^{-C_0 x} x \, dx$$

$$\leq \frac{2\pi}{k^2} C_1 e^{-C_0 k}.$$

From Markov's inequality, $\mathbb{P}(I_k > 0) = \mathbb{P}(I_k \geq 1/k^2) \leq k^2 \mathbb{E}(I_k)$ and from the Borel–Cantelli lemma, we obtain that almost surely $I_k$ is equal to 0 for $k$ large enough.

Similarly, let $p'$ be such that $q = p(1 - \varepsilon) < p' < p$, chosen as in the proof of Theorem 5.9. Using equations (43), (44) and (45) we get similarly

$$
\begin{aligned}
\mathbb{E}(L_k) &= \frac{2\pi}{k^2} \int_0^{(1-\varepsilon)pk} \mathbb{P}(H(xe_x) \geq k) x\, dx \\
&\leq \frac{2\pi}{k^2} \left( \int_0^{(1-\varepsilon)pk} x C_1 e^{-C_0(kp'-x)} x\, dx + C_1 e^{-C_0 k} \int_0^{(1-\varepsilon)pk} x\, dx \right) \\
&\leq \frac{2\pi}{k^2} C_1 e^{-C_0 k}.
\end{aligned}
$$

Using Markov's inequality and the Borel–Cantelli lemma, we deduce that, almost surely, $L_k$ is 0 for $k$ large enough.

The ergodic properties of the PPP imply that

$$
J_k = k^{-2} N(B(O, (1+\varepsilon)pk) \setminus B(O, (1-\varepsilon)pk))
$$

converges almost surely and in mean toward $4\pi p\varepsilon$. Notice that $N(B(O, (1+\varepsilon)pk) \setminus B(O, (1-\varepsilon)pk))$ is not an increasing sequence of convex sets. To prove this convergence, we need to use the independent properties of the PPP.

We have proved that, for all $\varepsilon > 0$, almost surely,

$$
\limsup_k \frac{|G_k - N(B(O, pk))|}{k^2} \leq 4\pi p\varepsilon.
$$

Hence, almost surely,

$$
\lim_k \frac{G_k}{k^2} = \lim_k \frac{N(B(O, p))}{k^2} = \pi p^2.
$$

The proof for the $L^1$ convergence is a consequence of the dominated convergence theorem.

Equation (1) holds since we have seen that a.s. for $k$ large enough $I_k$ and $L_k$ are both equal to 0, where $I_k$ is the cardinal of $\mathcal{T}(k) \cap B(O, p+\varepsilon)^c$ and $L_k$ is the cardinal of $\mathcal{T}(k)^c \cap B(O, p-\varepsilon)$.

**6. Spatial averages of edge lengths.** Consider the total edge length of the RST for points included in the ball $B(O, x)$:

$$
\mathcal{L}_x := \sum_{X \in N} \mathbb{1}(X \in B(O, x)) |X - \mathcal{A}(X)|.
$$

It is well known that for the minimal spanning tree, the subadditive ergodic theorem allows one to prove that $\frac{\mathcal{L}_x}{x^2}$ tends almost surely toward a constant (see [26]). We prove that the same holds for the RST (with, of course, a larger constant).

From Campbell's formula,

$$\mathbb{E}\mathcal{L}_x = 2\pi \int_0^x \mathbb{E}(L(t))t\,dt.$$

With the change of variable $u = \frac{t}{x}$, this leads to

$$\mathbb{E}\frac{\mathcal{L}_x}{x^2} = 2\pi \int_0^1 u\mathbb{E}(L(xu))\,du.$$

The dominated convergence theorem together with equation (10) gives

$$\lim_{x\to\infty} \mathbb{E}\frac{\mathcal{L}_x}{x^2} = 2\pi \int_0^1 u\frac{1}{\sqrt{2}}\,du = \pi/\sqrt{2}.$$

We will prove a stronger result: $\frac{\mathcal{L}_x}{x^2}$ converges almost surely and in $L^1$ toward $\pi/\sqrt{2}$.

To prove this, we consider here a slightly different problem. We sample $n$ points uniformly and independently on the unit disk. This defines a finite point set $F_n = \{O, X_1, \ldots, X_n\}$. We can then construct the RST associated to this point set. The total edge length of this RST is

$$\mathcal{L}(F_N) := \sum_{k=1}^N |X_k - \mathcal{A}(X_k)|.$$

First notice that $\mathcal{L}$ is homogeneous of order 1: for all sets $F_n$ as above and all positive real numbers $r$, we have $\mathcal{L}(rF_n) = r\mathcal{L}(F_n)$, where $rF_n = \{O, rX_1, \ldots, rX_n\}$. Using this and the fact that the ratio $N(B(O,x))/x^2$ tends a.s. toward $\pi$, it is easy to check that $\frac{\mathcal{L}_x}{x^2}$ converges a.s. toward $\pi/\sqrt{2}$ if and only if $\frac{\mathcal{L}(F_n)}{\sqrt{n}}$ tends to $\sqrt{\frac{\pi}{2}}$. Since we have already computed the mean of $\frac{\mathcal{L}_x}{x^2}$, it is sufficient to prove that $\frac{\mathcal{L}(F_n)}{\sqrt{n}}$ converges a.s. toward a constant.

For proving the last property, we use the smoothness of $\mathcal{L}(F_n)$ (Theorem 6.1) and the Rhee and Talagrand concentration inequalities.

THEOREM 6.1. *There exists a positive constant $C_1$ such that, for all finite subsets $F$ and $G$ as defined above,*

$$|\mathcal{L}(F \cup G) - \mathcal{L}(F)| \le C_1\sqrt{|G|},$$

*in particular, $\mathcal{L}(F) \le C_1\sqrt{|F|}$.*

PROOF.   $\mathcal{L}$ clearly satisfies the subadditive property: for all finite subsets $F$ and $G$,

(47)                    $$\mathcal{L}(F \cup G) \le \mathcal{L}(G) + \mathcal{L}(F).$$

From Lemma 3.4.1 of [27], we deduce that there exists a constant $C_1$ such that $\mathcal{L}(F) \leq C_1\sqrt{|F|}$. Subadditivity then implies

$$\mathcal{L}(F) \geq \mathcal{L}(F \cup G) - \mathcal{L}(G)$$
$$\geq \mathcal{L}(F \cup G) - C_1\sqrt{|G|}.$$

It remains to prove that $\mathcal{L}(F) \leq \mathcal{L}(F \cup G) + C_1\sqrt{|G|}$, for all finite sets $F$ and $G$ as above.

Let $Y \in G$ and suppose that the points $X_1, \ldots, X_n$ of $F$ all have $Y$ as an ancestor in the RST built over the set $\{O\} \cup F \cup G$. In particular, $|X_i| \geq |Y|$. Suppose $|X_k| \geq |X_j|$, then

$$|X_k - Y|^2 \leq |X_k - X_j|^2$$
$$= |X_k - Y|^2 + |X_j - Y|^2 - 2|X_k - Y||X_j - Y|\cos\widehat{X_kYX_j}.$$

Thus, if $\widehat{X_kYX_j} \leq \frac{\pi}{3}$, $|X_k - Y| \leq |X_j - Y|$. The inequality $|X_k - Y|^2 + |X_j - Y|^2 - 2|X_k - Y||X_j - Y|\cos\widehat{X_kYX_j} \leq |X_j - Y|^2$ holds for $|X_k - Y| \in [0, 2|X_j - Y|\cos\widehat{X_kYX_j}] \supset [0, |X_j - Y|]$. It follows that

$$|X_k| \geq |X_j| \quad \text{and} \quad \widehat{X_kYX_j} \leq \frac{\pi}{3}$$

(48) $$\text{imply}$$

$$|X_k - Y| \leq |X_j - Y| \quad \text{and} \quad |X_k - X_j| \leq |X_j - Y|.$$

Let $\Theta(Y, X)$ denote the oriented angle between $\vec{OY}$ and $\vec{OX}$. Due to the origin, if $X$ is connected to $Y \neq O$ in the RST, then $|\Theta(Y, X)| \leq \pi/2$.

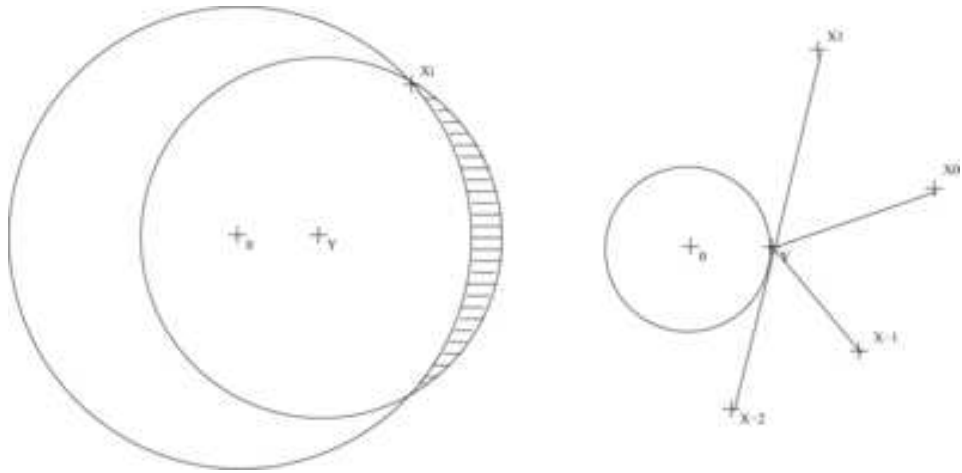

Fig. 12. *On the* left, *the dashed area is* $\mathcal{G}(X_i)$, *on the* right, *the set of connected points to* $Y$.

We index the $n$ points of $F$ connected to $Y$ by their increasing oriented angle $\Theta(Y, X)$ such that $|\Theta(Y, X_0)|$ is minimal. $X_1, \ldots, X_d$ are the points counted clockwise from $X_0$: $\Theta(Y, X_0) \leq \Theta(Y, X_1) \leq \cdots \leq \Theta(Y, X_1) \leq \pi/2$ and $X_{-1}, \ldots, X_{-d'}$ are the points counted counter-clockwise from $X_0$ with $n = d + d' + 1$ (see right picture in Figure 12). We need a tie-breaking rule: if $\Theta(Y, X) = \Theta(Y, X')$, the point with the higher norm has an index closer to $O$.

Assume now that, for a given $j \geq 0$, we have both $|X_{j+1}| \geq |X_j|$ and $\widehat{X_j Y X_{j+1}} \leq \frac{\pi}{3}$. Then from equation (48), $|X_{j+1} - Y| \leq |X_j - Y|$, $X_{j+1}$ belongs to the set $\mathcal{G}(X_j)$ of points closer to $Y$ than $X_j$ and with a norm larger than $X_j$. However, by elementary considerations, $X \in \mathcal{G}(X_j)$ implies that $|\Theta(Y, X)| \leq |\Theta(Y, X_j)|$ (see left picture in Figure 12) and this contradicts $\Theta(Y, X_{j+1}) > \Theta(Y, X_j)$ (the strict inequality comes from the tie-breaking rule). Similarly, for $j \leq 0$, if $\widehat{X_{j-1} Y X_j} \leq \frac{\pi}{3}$, then $|X_{j-1}| \leq |X_j|$.

There are at most 6 points such that $\widehat{X_i Y X_{i+1}} \geq \frac{\pi}{3}$ $(i \geq 0)$ or $\widehat{X_i Y X_{i-1}} \geq \frac{\pi}{3}$ $(i \leq 0)$. Therefore, from equation (48), there are at most 6 points such that $(i \geq 0$ and $|X_{i+1}| \leq |X_i|$ and $|X_{i-1} - X_i| \leq |X_i - Y|)$ or $(i < 0$ and $|X_{i+1}| \geq |X_i|$ and $|X_{i+1} - X_i| \leq |X_i - Y|)$ does not hold.

Let $\mathcal{A}_F(X)$ denote the ancestor of $X$ in the RST built on the set $\{O\} \cup F$. We have

$$\begin{aligned}
\mathcal{L}(F) &= \sum_{X \in F} |X - \mathcal{A}_F(X)| \\
&= \sum_{X \in F} \mathbb{1}(\mathcal{A}_F(X) = \mathcal{A}_{F \cup G}(X)) |X - \mathcal{A}_{F \cup G}(X)| \\
&\quad + \sum_{Y \in G} \sum_{X \in F \cap \mathcal{A}_{F \cup G}^{-1}(Y)} |X - \mathcal{A}_F(X)| \\
&\leq \mathcal{L}(F \cup G) + \sum_{Y \in G} \sum_{X \in F \cap \mathcal{A}_{F \cup G}^{-1}(Y)} |X - \mathcal{A}_F(X)| - |X - Y|.
\end{aligned}$$

If $F \cap \mathcal{A}_{F \cup G}^{-1}(Y) = \{X_{-d'}, \ldots, X_0, \ldots, X_d\}$, we have seen above that for at most 6 points $|X_i - \mathcal{A}_{F \cup G}(X_i)| = |X_i - Y| \geq |X_i - X_{i-1}| \geq |X_{i-1} - \mathcal{A}_F(X_{i-1})|$ $(i > 0)$ or $|X_i - Y| \geq |X_{i+1} - X_i| \geq |X_{i+1} - \mathcal{A}_F(X_{i+1})|$ $(i < 0)$ does not hold.

Henceforth, if $J$ denotes the set of points such that the preceding inequality does not hold and $H(Y) = \{X_d, X_{-d'}\} \cup J$, we have $|H(Y)| \leq 8$ and

$$\sum_{X \in F \cap \mathcal{A}_{F \cup G}^{-1}(Y) \setminus H(Y)} |X - \mathcal{A}_F(X)|$$

$$\leq \sum_{i=1}^{d} |X_i - X_{i-1}| + \sum_{i=2}^{d'} |X_{-i} - X_{-i+1}|$$

$$\leq \sum_{X \in F \cap \mathcal{A}_{F \cup G}^{-1}(Y)} |X - \mathcal{A}_{F \cup G}(X)|.$$

If $H = \bigcup_{Y \in G} H(Y)$, we have $|H| \leq 8|G|$. Using subadditivity, we deduce

$$\mathcal{L}(F) \leq \mathcal{L}(F \cup G) + \mathcal{L}(H) \leq \mathcal{L}(F \cup G) + C\sqrt{8|G|}. \qquad \square$$

Theorem 6.1 ensures that we can apply the Rhee and Talagrand concentration inequalities to the functional $\mathcal{L}$ (Theorem 1 of [25] and Theorem 11.3.2 of [28]): $\frac{\mathcal{L}(F_N)}{\sqrt{N}}$ converges a.s. toward its mean. Finally, we have proved that a.s. and in $L^1$,

$$(49) \qquad \lim_{x \to \infty} \frac{\mathcal{L}_x}{x^2} = \pi/\sqrt{2}.$$

More generally, we may also consider a power-weighted edge and for $\alpha \geq 0$, define

$$\mathcal{L}^{(\alpha)}(F) := \sum_{X \in F} |X - \mathcal{A}(X)|^\alpha.$$

The proof of Theorem 6.1 is unchanged if we replace $\mathcal{L}^{(\alpha)}$ by $\mathcal{L}$ (however, the constant $C_1$ does depend on $\alpha$). Define

$$(50) \qquad \lambda_\alpha = \alpha \int_0^\infty r^{\alpha-1} e^{-\pi r^2/2} \, dr = \mathbb{E}|O - \mathcal{A}_{e_1}(O)|^\alpha.$$

From equation (10), using the Campbell formula, we get that

$$\lim_{x \to +\infty} \mathbb{E}\mathcal{L}^{(\alpha)}(x)/x^2 = \pi \lambda_\alpha.$$

We finally deduce the following:

THEOREM 6.2. *For all $\alpha \geq 0$ a.s. and in $L^1$,*

$$\lim_{x \to \infty} \frac{\mathcal{L}_x^{(\alpha)}}{x^2} = \pi \lambda_\alpha.$$

We can rewrite this result as $1/N(B(O, x)) \sum_{X \in N} |X - \mathcal{A}(X)|^\alpha$ tends a.s. toward $\lambda_\alpha$. That is, the spatial average of the power-weighted lengths of the edges tends toward the distribution of the length of $(O, \mathcal{A}_{e_1}(O))$ in the DSF.

It is crucial to notice that $\lambda_\alpha < l_\alpha$: the spatial average and the average along a long path do not coincide.

REMARK 6.3. Theorem 2.1 of [24] gives the weak convergence of $\mathcal{L}_x^{(\alpha)}/x^2$ to $\pi \lambda_\alpha$. We have done our analysis on the length of an edge, the same type of result could be obtained for other stabilizing functionals. In order to derive weak laws for stabilizing functionals, we can directly invoke Theorem 2.1 of [24].

## 7. Model extension and open problems.

7.1. *Greedy forests.* The radial spanning tree lies in a large class of spanning forests which are locally defined. We could extend the definition of the radial spanning tree over a point set $N$ in $\mathbb{R}^d$ as follows.

Let $l$ be a measurable function from $\mathbb{R}^d$ to $\mathbb{R}_+$ and $L$ be a measurable function from $\mathbb{R}^d \times \mathbb{R}^d$ to $\mathbb{R}_+$. Suppose that for all $x, y, z, t \in N$, $\{X, Y\} \neq \{Z, T\}$, $L(X, Y) \neq L(Z, T)$, for $X \neq Y$, $l(X) \neq l(Y)$, and for all compact $K \subset \mathbb{R}_+$, $l^{-1}(K) \cap N$ is finite. Then we can define the following forest $\mathcal{F} = (N, E)$: for $l(Y) < l(X)$, $(X, Y) \in E$ if and only if $Y = \arg\min_{Z \in N, l(Z) < l(X)} L(X, Z)$.

When $l(X) = |X|$, $L(X, Y) = |X - Y|$ and $O \in N$, we define the radial spanning tree; if $l(X) = \langle X, e_x \rangle$ and $L(X, Y) = |X - Y|$, this is the directed spanning forest; if instead $L(X, Y)$ is the length of a well-chosen cylinder, we obtain the Poisson forest of [12].

Let $R(X)$ denotes the set of successive ancestors of vertex at $X$. The construction of $R(X)$ can be thought of as greedy because at each vertex of $R(X)$ a one-step optimization is performed.

As an example, in the Appendix, Section A.2, we discuss briefly the radial spanning tree obtained by choosing the $L^\infty$-norm.

7.2. *The radial spanning tree of a Voronoi cell.* An interesting way to extend the RST is to consider two independent Poisson point processes, $N_0 = \{T_n^0\}$, the point process of *cluster heads* (we use here the terminology of sensor networks, which motivate this extension) of intensity $\lambda_0$ and $N_1 = \{T_n^1\}$, the point process of *nodes*, with intensity $\lambda_1$. The first point process tessellates the plane in Voronoi cells. We denote by $V_n$ the Voronoi cell of point $T_n^0$ w.r.t. the points of $N_0$. Two forests can then be defined in relation with this tessellation:

- The family of *internal* RSTs: the $n$th tree of this forest, $\mathcal{T}_n$, is the RST built using the points of $N_1$ that are contained in $V_n$, with $T_n^0$ as a root.
- The family of *local* RSTs: if node $X$ belongs to $V_n$, one defines its ancestor as the point of $(N_1 \cup \{T_n^0\}) \cap B(T_n^0, |X - T_n^0|)$ that is the closest to $X$. Notice that this ancestor does not necessarily belong to $V_n$. Nevertheless, this rule defines a forest too (see Lemma 7.1). One then defines the $n$th local RST tree $\mathcal{U}_n$ as the tree which is the union of all the paths from nodes with ultimate ancestor $T_n^0$.

In what follows we concentrate on the second case which can be analyzed using the same type of tools as in Section 3. Figure 13 depicts a sample of such a forest.

LEMMA 7.1. *Almost surely, there exists no node $X$ of $N_0$ such that the sequence of ancestors of $X$ based on the local RST rule contains node $X$.*

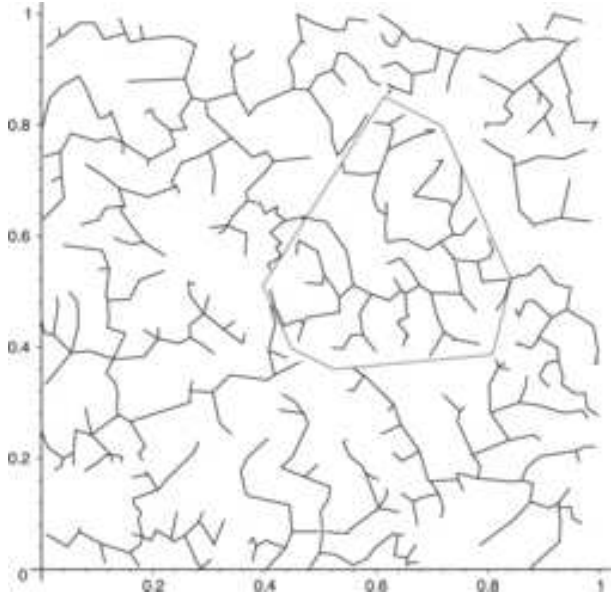

Fig. 13. *Local Voronoi radial spanning trees of* 600 *nodes uniformly and independently distributed in the unit square, w.r.t.* 10 *cluster heads, also uniformly and independently distributed in the unit square. The Voronoi cell of one of the cluster heads is depicted in red.*

PROOF. Let $\{Y_i\}_{i\geq 0}$ be the sequence of ancestors of $X = Y_0$. If $Y_i = X$ for some $i > 0$, then necessarily the points of $\{Y_i\}$ belong to different Voronoi cells (if this were not the case, then the distance to the cluster head of the cell to which all points belong would be strictly decreasing, which forbids cycles). Then one can then rewrite $\{Y_i\}_{i\geq 0}$ as

$$\{Y_i\}_{i\geq 0} = \{Z_0(1), \dots, Z_0(n_0), Z_1(1), \dots, Z_1(n_1), \dots, Z_j(1), \dots, Z_j(k), \dots\},$$

with $Z_l(1), \dots, Z_l(n_l) \in W_l$ for all $0 \leq l \leq j$, where $\{W_j\}_{j\geq 0}$ is a sequence of cells such that $W_l \neq W_{l+1}$ for all $l < j$, $n_l$ is a sequence of integers and $Z_j(k) = X$. Let $S_l$ denote the cluster head of $W_l$. Then the definition of the local RST implies that a.s.

$$|X - S_0| = |Z_0(1) - S_0| > |Z_0(2) - S_0| > \cdots > |Z_0(n_0) - S_0| > |Z_1(1) - S_0|.$$

Since $Z_1(1)$ belongs to $W_1$, we have a.s.

$$|Z_1(1) - S_0| > |Z_1(1) - S_1|.$$

For the same reasons, for all $l = 1, \dots, j-1$,

$$|Z_l(1) - S_l| > |Z_l(2) - S_l| > \cdots > |Z_l(n_l) - S_l| > |Z_{l+1}(1) - S_l|$$

and

$$|Z_{l+1}(1) - S_l| > |Z_{l+1}(1) - S_{l+1}|, \qquad \text{a.s.}$$

In addition,

$$|Z_j(1) - S_j| > |Z_l(2) - S_l| > \cdots > |Z_j(k) - S_j| = |X - S_0|.$$

Hence, a contradiction. $\quad\square$

Let $\mathcal{L}_n$ denote the total length of all edges from vertices in $V_n$. Let $\mathbb{E}_0$ denote the Palm probability w.r.t. $N_0$. We have

$$\mathcal{L}_n = \sum_m \mathbb{1}(T_m^1 \in V_n)L_m,$$

where $L_m$ is the length of the link that connects $T_m^1$ to its ancestor. Using the fact that $T_m^1 \in V_n$ iff $N_0(B(T_m^1, |T_m^1 - T_n^0|)) = 0$ and the fact that $L_m > u$ with $u < |T_m^1 - T_n^0|$ iff $N_1(B(T_m^1, u) \cap B(T_n^0, |T_m^1 - T_n^0|)) = 0$, we get from Campbell's formula that

$$\mathbb{E}_0\left(\sum_m \mathbb{1}(T_m^1 \in V_n)L_m\right) = 2\pi\lambda_1 \int_{r=0}^\infty e^{-\lambda_0\pi r^2}\left(\int_{u=0}^r e^{-\lambda_1 M(r,u)}\,du\right)r\,dr,$$

with $M(r,u)$ the lens defined in Section 3.1. Hence,

$$(51) \qquad \mathbb{E}_0(\mathcal{L}_0) = 2\pi\lambda_1 \int_{r=0}^\infty e^{-\lambda_0\pi r^2}\left(\int_{u=0}^r e^{-\lambda_1 M(r,u)}\,du\right)r\,dr.$$

7.3. *Open problems.* The local geometry of the RST is rather well understood. Unfortunately, the distribution of the degree of a vertex is still unknown. It would be appealing to compute this distribution at least in the DSF. In contrast with what happens in the minimal spanning tree, the degree is not upper bounded and so the moments of this distribution could be large.

Properly scaled, the path $R(X)$ of successive ancestors of $X$ in the DSF converges weakly toward the Brownian motion. An interesting problem is to find a functional central limit theorem for $R_O(X)$. Along this line, we may prove that the DSF converges weakly toward the Brownian web. Proving a weak limit for the RST is a challenging question.

## APPENDIX

**A.1. Distribution for the $n$th point of the PPP.** We will denote by $T_n$ the $n$th point of the Poisson point process when sorting points according to their distance to the origin. For a bounded Borel set $A$, let

$$N(A) = \sum_n \mathbb{1}(T_n \in A).$$

We will denote the angle of $T_n$ by $\phi_n$ and its norm by $|T_n| = \nu_n$.

The sequence $\{\phi_n\}$ is i.i.d., uniform on $(0, 2\pi)$ and independent of the sequence $\{\nu_n\}$.

The random variable $\nu_n^2$ has a Gamma distribution of parameter $(n, \pi\lambda)$ and the sequence $\{\nu_n\}$ is a Markov chain with kernel of density

$$K(y, t) = 2\pi\lambda t e^{-\lambda\pi(t^2 - y^2)}, \qquad t \geq y.$$

Indeed, the Markov property is immediate from the definition of Poisson point processes. In addition, by a direct martingale argument,

$$\mathbb{E}(e^{\lambda(1-z)\pi(\nu_{n+1}^2 - y^2)}|\nu_n = y) = z^{-1}.$$

Hence,

$$\mathbb{E}(e^{-u\nu_{n+1}^2}|\nu_n = y) = e^{-uy^2}\frac{\lambda\pi}{u + \lambda\pi},$$

so that conditioned on $\nu_n = y$, $\nu_{n+1}^2$ is the sum of $y^2$ and of an exponential random variable of parameter $\lambda\pi$. Hence, $\nu_n^2$ has a Gamma distribution of parameter $(n, \pi\lambda)$ indeed. In addition, the conditional density of $\nu_{n+1}$ given $\nu_n = y$ is

$$f_n(t, y) = 2\pi\lambda t e^{-\lambda\pi(t^2 - y^2)}, \qquad t \geq y.$$

Let $\gamma_{n,\pi}$ be the density of a Gamma distribution of parameter $(n, \pi)$; letting $L_n = |T_n - \mathcal{A}(T_n)|$, we have

$$\mathbb{P}(L_n \geq r) = \mathbb{P}(\nu_n \geq r; N(B(T_n, r) \cap B(O, \nu_n)) = 0)$$

$$= \int_{r^2}^{+\infty} \mathbb{P}(N(B(T_n, r) \cap B(O, \sqrt{t})) = 0|\nu_n^2 = t)\gamma_{n,\pi}(t)\, dt$$

$$= \int_{r^2}^{+\infty} \left(1 - \frac{M(\sqrt{t}, r)}{\pi t^2}\right)^{n-1} \gamma_{n,\pi}(t)\, dt.$$

The same type of computation can be done on the distribution of the edge $(T_n, \mathcal{A}(T_n))$ or for $D_n = \mathbb{E}D(T_n)$.

**A.2. Radial spanning tree with $L^\infty$-norm.** In the $L_\infty$ case, if $X = (t\cos\theta, t\sin\theta)$, then for $r \leq |X|_\infty = t\max(|\cos(\theta)|, |\sin(\theta)|)$,

$$M(X, r) = 2r^2\mathbb{1}(r < g(X)) + (r^2 + rg(X))\mathbb{1}(r \geq g(X))$$

with

$$g(X) = t(\max(|\cos\theta|, |\sin\theta|) - \min(|\cos\theta|, |\sin\theta|))$$

(see Figure 14).

In addition, when denoting by $L(X)$ the $L_\infty$ distance of point $X$ to its parent vertex in the tree, we still have

$$\mathbb{P}(L(X) \geq r) = \mathbb{1}(r \leq |X|_\infty)e^{-M(X,r)}.$$

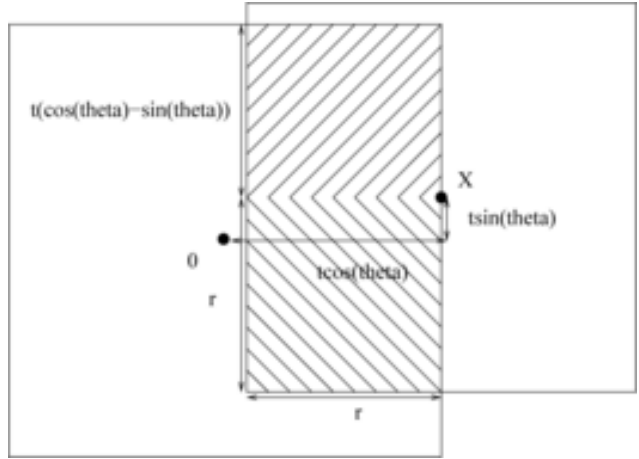

FIG. 14.    *The $L_\infty$ "lens."*

### A.3. Tail inequality for sums of nonnegative random variables.

LEMMA A.1.    *Let $(X_k)_k$ be a sequence of real random variables and $\{\mathcal{F}_k\}_k$ a filtration of this process. If $\mathbb{E}(X_{n+1}|\mathcal{F}_n) = 0$ and $\mathbb{P}(|X_{n+1}| \geq t|\mathcal{F}_n) \leq C_1 \exp(-C_0 t)$, then for all $t_0 > 0$, there exists positive constants $C_0$, $C_1$ such that, for all $t \geq t_0$,*

$$\mathbb{P}\left(\left|\frac{1}{n}\sum_{k=1}^{n} X_k\right| \geq t\right) \leq C_1 e^{-C_0 n t}.$$

PROOF.    This lemma relies on a classical computation on large deviations, we only give a sketch of the proof. Let $\Lambda_n(\lambda) = \ln \mathbb{E}(e^{\lambda X_n}|\mathcal{F}_{n-1})$ and $\Lambda_n^*(t) = \sup_\lambda \lambda t - \Lambda_n(\lambda)$, the Fenchel–Legendre transform of $\Lambda_n$. With $\Lambda_n^*(0) = 0$, the condition $\mathbb{P}(|X_{n+1}| \geq t|\mathcal{F}_n) \leq C_1 \exp(-C_0 t)$ ensures that, for $t \neq 0$, $\Lambda_n^*(t)$ is positive and lower bounded uniformly in $n$ by a positive $C_0$. For $t \geq 0$, $\Lambda_n^*$ is nondecreasing and convex (refer to [9]), hence, if $t \geq t_0 > 0$, $\Lambda_n^*(t) \geq t\frac{\Lambda_n^*(t_0)}{t_0} \geq tC_0$. As usual for upper bounds in large deviation, the rest of the proof follows from Chernoff's inequality.    □

### A.4. Some tail inequalities in the $GI/GI/\infty$ queue.    This section uses the classical notation of queueing theory. This leads to conflicts with some of the notation used elsewhere in the paper. This should not lead to difficulties, as only the results stated in Lemma A.2 will actually be used.

Let $\{\sigma_n, \tau_n\}, n \in \mathbb{Z}$, be an i.i.d. sequence of $\mathbb{R}_+ \times \mathbb{R}_+$-valued random variables representing the service times and inter-arrival times in a $GI/GI/\infty$ queue. For $n$ fixed, the random variables $\sigma_n$ and $\tau_n$ are possibly dependent.

We set $T_0 = 0$ as the arrival time of customer 0; for $n \geq 1$, $T_n = \sum_{k=0}^{n-1} \tau_k$ is the arrival time of the $n$th customer. Let $Y \in \mathbb{R}_+$ be a nonnegative initial condition, independent of the $\{\sigma_n, \tau_n\}$ sequence. We set $W_0^{[Y]} = Y$, and for $n \geq 1$, we define

$$W_n^{[Y]} = \max\left( Y - \sum_{k=0}^{n-1} \tau_k, \max_{2 \leq i \leq n} \sigma_{i-1} - \sum_{k=i-1}^{n-1} \tau_k \right)^+$$

$$= \max\left( Y - (T_n - T_0), \max_{2 \leq i \leq n} \sigma_{i-1} - (T_n - T_{i-1}) \right)^+.$$

The random variable $W_n^{[Y]}$ is the largest residual service time just before the arrival of the $n$th customer in the $GI/GI/\infty$ queue with initial condition $Y$.

The following additional assumptions are made:

(i) There exist constants $C_0, C_1$ such that $\mathbb{P}(\sigma_1 \geq t) \leq C_1 \exp(-C_0 t)$ and $\mathbb{P}(Y > t) \leq C_1 \exp(-C_0 t)$.

(ii) The support of $\sigma_1$ is $\mathbb{R}_+$.

(iii) $\mathbb{P}(\tau_1 > 0) > 0$.

The Loynes' sequence $\{M_n\}$ of this $GI/GI/\infty$ queue is defined by $M_0 = 0$ and

$$M_n = \max_{-n+1 \leq i \leq 0} \left( \sigma_{i-1} - \sum_{k=i-1}^{-1} \tau_k \right)^+, \qquad n \geq 1.$$

This sequence is nondecreasing in $n$ and it a.s. converges to

$$(52) \qquad M = \sup_{i \leq 0} \left( \sigma_{i-1} - \sum_{k=i-1}^{-1} \tau_k \right)^+.$$

The random variable $M$ is a.s. finite. Indeed, we can easily obtain a stronger assertion. Indeed, by assumption (i), for some $s > 0$,

$$\mathbb{E} \exp(sM) \leq 1 + \sum_{i \leq 0} \mathbb{E} \exp\left( s \left( \sigma_{i-1} - \sum_{k=i-1}^{-1} \tau_k \right) \right)$$

$$\leq 1 + \mathbb{E} \exp(s\sigma_1) \sum_{i \leq 0} \mathbb{E} \exp(-si\tau_1) < \infty.$$

From this and Chernoff's inequality, we deduce that

$$(53) \qquad \mathbb{P}(M > t) \leq C_1 e^{-C_0 t},$$

for some constants $C_0, C_1$.

Now, we define

$$(54) \qquad \nu(Y) = \nu = \inf\left\{n \geq 2 : Y - \sum_{k=0}^{n-1} \tau_k < 0\right\}.$$

From time $\nu$ on, the initial workload does not count anymore, that is, for $n \geq \nu$, $W_n^{[Y]} = \max_{2 \leq i \leq n}(\sigma_{i-1} - T_n + T_{i-1})^+$. Note that $\nu$ has the same distribution as

$$\nu' = \max\left\{n \leq -1 : Y - \sum_{k=n}^{-1} \tau_k < 0\right\}.$$

More generally,

$$\left(\sum_{k=0}^{\nu-1} \tau_k, \sum_{k=1}^{\nu-1} \tau_k, \ldots, \tau_{\nu-1}\right)$$

and

$$\left(\sum_{k=-\nu'}^{-1} \tau_k, \sum_{k=-\nu'+1}^{-1} \tau_k, \ldots, \tau_{-1}\right)$$

have the same law, which implies that $M_{\nu'}$ and $W_\nu^{[Y]}$ have the same distribution. Since $M_{\nu'} \leq M$, we have

$$(55) \qquad W_\nu^{[Y]} \leq_{st} M.$$

Note that this bound is uniform in $Y$.

From assumption (iii), we may find $c > 0$ and $\varepsilon > 0$ such that $\mathbb{P}(\tau_1 \geq c) \geq \varepsilon$. Let $B_k = c\mathbb{1}(\tau_k \geq c)$. Using the independence between $\tau$ and $Y$ and Hoeffding's inequality,

$$\mathbb{P}(\nu > n) \leq \mathbb{P}\left(\sum_{k=0}^{n-1} B_k < Y\right) \leq \mathbb{E}\exp\left(-\frac{(c\varepsilon n - Y)^2}{2c^2 n}\right).$$

By assumption (i), $Y$ is such that $\mathbb{P}(Y > t) \leq C_1 \exp(-C_0 t)$, for some positive constants $C_0$, $C_1$, hence,

$$(56) \quad \mathbb{P}(\nu > n) \leq \mathbb{P}(Y > nt_0) + \exp\left(-\frac{(c\varepsilon n - nt_0)^2}{2c^2 n}\right) \leq C_1 \exp(-C_0 n),$$

for some positive constants $C_0, C_1$, uniformly on the initial conditions $Y$.

Let $\mathcal{F}_n$ be the $\sigma$-field generated by the random variables $Y$ and $\{(\sigma_k, \tau_k), k = 0, \ldots, n-1\}$. The sequence $\{W_n^{[Y]}\}$ is a $\{\mathcal{F}_n\}$-Markov chain and the random variables

$$\nu_{n+1} = \nu_n + \nu(W_{\nu_n}^{[Y]}), n \geq 1,$$

with $\nu(W)$ defined in (54) and with $\nu_1 = \nu = \nu(Y)$, are $\{\mathcal{F}_n\}$-stopping times. Using what precedes, one gets by induction that each $\nu_n$ is a.s. finite and that for all $n$,

$$(57) \quad P(\nu_{n+1} - \nu_n > m | \mathcal{F}_{\nu_n}) = P(\nu_{n+1} - \nu_n > m | W_{\nu_n}^{[Y]}) \leq C_1 e^{-C_0 m} \qquad \forall m,$$

$$(58) \qquad P(W_{\nu_{n+1}}^{[Y]} > x | \mathcal{F}_{\nu_n}) = P(W_{\nu_{n+1}}^{[Y]} > x | W_{\nu_n}^{[Y]}) \leq P(M > x) \qquad \forall x.$$

Using (57) and a Chernoff type bound, one gets

$$(59) \qquad\qquad P(\nu_n > \alpha n) \leq C_1 e^{-C_0 n},$$

for some positive constants $\alpha, C_0, C_1$.

Consider now the $\{\mathcal{F}_n\}$-stopping time

$$(60) \qquad \theta(Y) = \theta = \inf\{n \geq 1 : (\sigma_{n-1} - \tau_{n-1})^+ = W_n^{[Y]}\}.$$

Since $W_\nu^{[Y]} \leq_{st} M$, by assumption (ii) implies

$$P(\theta \leq \nu_1 + 2) \leq P(M + \sigma_{\nu_1} - \tau_{\nu_1} - \tau_{\nu_1+1} \leq \sigma_{\nu_1+1} - \tau_{\nu_1+1}) = \delta > 0.$$

In the same vein, when using (59) and (58), one gets that

$$P(\theta > n) \leq (1 - \delta)^{\lfloor n/\alpha \rfloor} + P(\nu_{\lfloor n/\alpha \rfloor} > n) \leq C_1 e^{-C_0 n},$$

for some positive constants $C_0, C_1$.

We define the sequence $\theta_0 = 0$, $\theta_1 = \theta(Y)$ and

$$(61) \qquad \theta_{n+1} = \inf\{k > \theta_n : \sigma_{k-1} - \tau_{k-1} = W_k^{[Y]}\}, \qquad n \geq 1.$$

Let $\xi_n = \mathbb{1}(W_{\nu_n+2}^{[Y]} = \sigma_{\nu_n+1} - \tau_{\nu_n+1}))^+$. We have $P(\xi_n = 1 | \mathcal{F}_{\nu_{n-1}}) \geq \delta$ for all $n$. Letting $0 < \varepsilon < \delta$, we have

$$(62) \begin{aligned} P(\theta_n > (\delta - \varepsilon)n/\alpha) &\leq P(\nu_{\lfloor n/\alpha \rfloor} > n) + P\left(\sum_{k=1}^{\lfloor n/\alpha \rfloor} \xi_k < n(\delta - \varepsilon)/\alpha\right) \\ &\leq C_1 e^{-C_0 n}, \end{aligned}$$

when using (59) and a Chernoff type bound.

We have proved the following lemma:

LEMMA A.2. *Let $\theta_n$ be the stopping time defined in equation* (61). *Under the foregoing probabilistic assumptions on $Y$, $(\tau_n)$ and $(\sigma_n)$, there exists $\beta > 0$ such that, for all $k$,*

$$P(\theta_{k+1} - \theta_k > n) \leq C_1 e^{-C_0 n} \quad\text{and}\quad P(\theta_n > \beta n) \leq C_1 e^{-C_0 n},$$

*for some positive constants $C_0, C_1$.*

## MATHEMATICAL NOTATION

- $\mathcal{A}(X)$: the ancestor of point at $X$ in the RST—see Section 2.2;
- $\mathcal{A}^k(X)$: the ancestor of level $k$ of point at $X$ in the tree—see Section 2.2;
- $\mathcal{A}_{e_1}(X)$: the ancestor of point at $X$ in the DSF with direction $e_1$;
- $B(x,r)$: the open ball of radius $r$ and center $x$;
- $C(X)$: the number of links of the tree that cross the circle of radius $|X|$—see Section 3.5;
- $C_0$: a small constant—see Section 2.2;
- $C_1$: a large constant—see Section 2.2;
- $D(X)$: the degree of a point located at $X$—see Section 2.3.2;
- $D_n$: the domain (13) associated with the $n$th point of the path;
- $\Delta(x)$: the maximal deviation of $R_0(x)$ and the $x$ axis—see equation (31);
- $\Delta(x,x')$: the maximal deviation of $R(x)$ between $x'$ and $x$ defined in equation (30);
- $\{\Phi_n\}$: Markov chain with state in the space of finite point processes—equation (17);
- $\mathcal{H}$: the right half-plane ($x>0$);
- $\mathcal{H}^-(t)$: the left half-plane ($x<t$);
- $H(X)$: the number of hops from point $X$ to the origin—see Section 2.3.2;
- $L(X)$: the length of the link connecting a point located at $X$ to its ancestor in the RST—see Section 2.3.2;
- $N$: the underlying Poisson point process;
- $p$: the mean stationary progress in a long path of the DSF—equation (27);
- $P(X)$: the progress of the link from a point at $X$ toward the origin—see Section 2.3.2;
- $R(x)$: the path from $(x,0)$ in the DSF—see Section 4.3;
- $R_0(X)$: the path from $X$ to $0$ in the RST—see Section 2.2;
- $\mathcal{S}$: the small set of the Markov chain $\{\Phi_n\}$—equation (22);
- $\sigma_A$: return time of the Markov chain $\{\Phi_n\}$ to set $A$—equation (18);
- $\mathcal{T}$: the RST;
- $\mathcal{T}_{e_1}$: the DSF with direction $e_1$;
- $\mathcal{T}(k)$: the subtree of the RST with vertices less than $k$ generations away from the root;
- $\{T_n\}$: the points of $N$;
- $\tau_n$: the Markovian times in the sequence of successive ancestors of a point—equation (15);
- $\xi_n$: the projection of $\mathcal{H} \setminus D_n$ on the $x$ axis—equation (16).

**Acknowledgments.** The authors would like to thank Marc Lelarge and Venkat Anantharam for their comments and suggestions on this paper.

Département d'Informatique
École Normale Supérieure
45 rue d'Ulm
F-75230 Paris Cedex 05
France
E-mail: francois.baccelli@ens.fr
          charles.bordenave@ens.fr